\documentstyle{amltd2004}
\begin{document}
\annalsline{155}{2002}
\received{May 13, 1999}
\revised{March 28, 2001}
\startingpage{317}
\def\bye{\end{document}}
 \font\tenrm=cmr10
\def\et{\endproclaim}
\def\ec{\endproclaim}
\def\el{\endproclaim}
\def\er{\enddemo}
\def\lrc#1{\left\{#1\right\}}
\def\lrp#1{\left(#1\right)}
\def\lrb#1{\left[#1\right]}
\catcode`\@=11
\font\twelvemsb=msbm10 scaled 1100
\font\tenmsb=msbm10
\font\ninemsb=msbm10 scaled 800
\newfam\msbfam
\textfont\msbfam=\twelvemsb  \scriptfont\msbfam=\ninemsb
  \scriptscriptfont\msbfam=\ninemsb
\def\msb@{\hexnumber@\msbfam}
\def\Bbb{\relax\ifmmode\let\next\Bbb@\else
 \def\next{\errmessage{Use \string\Bbb\space only in math
mode}}\fi\next}
\def\Bbb@#1{{\Bbb@@{#1}}}
\def\Bbb@@#1{\fam\msbfam#1}
\catcode`\@=12

 \catcode`\@=11
\font\twelveeuf=eufm10 scaled 1100
\font\teneuf=eufm10
\font\nineeuf=eufm7 scaled 1100
\newfam\euffam
\textfont\euffam=\twelveeuf  \scriptfont\euffam=\teneuf
  \scriptscriptfont\euffam=\nineeuf
\def\euf@{\hexnumber@\euffam}
\def\frak{\relax\ifmmode\let\next\frak@\else
 \def\next{\errmessage{Use \string\frak\space only in math
mode}}\fi\next}
\def\frak@#1{{\frak@@{#1}}}
\def\frak@@#1{\fam\euffam#1}
\catcode`\@=12


\def\ee{\begin{equation}}
\def\aa{\begin{eqnarray}}
\def\y{\begin{eqnarray*}}
\def\bd{\begin{description}}
 \def\eee{\end{equation}}
\def\eaa{\end{eqnarray}}
\def\ey{\end{eqnarray*}}
\def\ebd{\end{description}}
\def\nn{\nonumber}
 \def\<{\langle}
\def\>{\rangle}
\def\dm{\diamond}
\def\ox{\mbox{}}
\def\lb{\label}
\def\bs{\setminus}
\def\R{{\bf R}}
\def\C{{\bf C}}
\def\Z{{\bf Z}}
\def\N{{\bf N}}
\def\U{{\bf U}}
\def\Q{{\bf Q}}
\def\T{{\bf T}}
\def\ga{{\gamma}}
\def\th{{\theta}}
\def\om{{\omega}}
\def\Om{{\Omega}}
\def\ep{{\varepsilon}}
\def\lm{{\lambda}}
\def\dl{{\delta}}
\def\sg{{\sigma}}
\def\Sg{{\Sigma}}
\def\vf{{\varphi}}
\def\F{{\cal F}}
\def\V{{\cal V}}
\def\G{{\cal G}}
\def\K{{\cal K}}
\def\P{{\cal P}}
\def\J{{\cal J}}
\def\Li{{\cal L}}
\def\I{{\cal I}}
\def\E{{\cal E}}
\def\M{{\cal M}}
\def\H{{\cal H}}
\def\Na{{\cal N}}
\def\ol#1{\overline{#1}}  
\def\td#1{\tilde{#1}}
\def\gl{{\rm gl}}
\def\GL{{\rm GL}}
\def\Sp{{\rm Sp}}
\def\mod{\;{\rm mod}\;}
\def\diag{{\rm diag}}
\def\ind{{\rm ind}}

\font\emi= cmmi10 scaled 1700 
\font\eightmi=cmmi10
\font\titr=cmr10 scaled 1700
\font\tmr=cmr10
\font\emr=cmr6

\title{Closed characteristics on compact\\ convex hypersurfaces
in  R\raise6pt\hbox{{\tmr \char50}{\eightmi \char110}}} 
\def\titleheadline#1{\def\one{#1}\ifx\one\empty\else
\gdef\thetitle{{\frenchspacing%
\let\\ \relax
{#1}}}\fi}
\newif\ifshort
\def\shortname#1{\global\shorttrue\xdef
\theauthors{{\eightsc\uppercase{#1}}}}
\let\shorttitle\titleheadline
\shorttitle{ \eightsc\uppercase{Compact convex hypersurfaces
in} {\eightpoint \bf R}\raise3pt\hbox{{\emr \char50}{\fivei \char110}}} 

\shorttitle{}

\acknowledgement{The first author was partially supported by the 973 Program of STM,
NNSF, MCME, RFDP, PMC Key Lab of EM of China, S. S. Chern Foundation,
Hong Kong Qiu Shi Sci. Tech. Foundation, and CEC of Tianjin.  The first author is an associate member of ICTP.  The
second author was partially supported by Hong Kong Qiu
Shi Sci.\ Tech.\ Foundation.}
  \twoauthors{Yiming Long}{Chaofeng Zhu}
 \institutions{Nankai Institute of Mathematics, Nankai University, Tianjin 300071, The People's Republic of
China\\ {\eightpoint {\it E-mail addresses\/}: longym@nankai.edu.cn, zhucf@nankai.edu.cn}}

\bigbreak \centerline{\bf Abstract}
\vglue12pt

 For any given compact $C^2$ hypersurface $\Sigma$ in $\R^{2n}$
bounding a strictly convex set with nonempty interior, in this paper
an invariant $\varrho_n(\Sigma)$ is defined and satisfies
$\varrho_n(\Sigma)\ge [n/2]+1$, where $[a]$ denotes the greatest
integer which is not greater than $a\in\R$. The following results
are proved in this paper. There always exist at least $\varrho_n(\Sigma)$
geometrically distinct closed characteristics on $\Sigma$. If all the
geometrically distinct closed characteristics on $\Sigma$ are
nondegenerate, then $\varrho_n(\Sigma)\ge n$. If the total
number of geometrically distinct closed characteristics on $\Sigma$ is
finite, there exists at least an elliptic one among them, and there exist at
least $\varrho_n(\Sigma)-1$ of them possessing irrational mean indices. If
this total number is at most $2\varrho_n(\Sigma)-2$, there exist at
least two elliptic ones among them. \vglue12pt

\section{Introduction and main results}\lb{s1}

1.1. {\it Main results}. 
 Let $\Sigma$ be a $C^2$-compact hypersurface in $\R^{2n}$ bounding
a strictly convex compact set $C$ with nonempty interior, where  $\Sigma$ has a
non-vanishing Gaussian curvature. In this paper we study  closed characteristics
on such hypersurfaces.  Without loss of generality, we may assume
$0\in C$. We denote the set of all such hypersurfaces in $\R^{2n}$ by
$\H(2n)$. For $x\in\Sigma$, let $N_{\Sigma}(x)$ be the outward normal unit
vector at $x$ of $\Sigma$. We consider the dynamics problem of finding
$\tau>0$ and an absolutely continuous curve $x\colon[0,\tau]\to\R^{2n}$
such that
\ee \left\{ \begin{array}{ll}
\dot{x}(t) &\hskip-8pt = JN_{\Sigma}(x(t)), \quad x(t)\in\Sigma,
                   \qquad\hbox{for all } t\in\R,\\
             x(\tau) &\hskip-8pt  = x(0),  \end{array}\right.
    \lb{eq1.1}\eee
where $J=\pmatrix{0&-I\cr I& 0\cr}$ is the standard
symplectic matrix on $\R^{2n}$. A solution $(\tau,x)$
of the problem (\ref{eq1.1}) is called a {\it closed characteristic}
on $\Sigma$. Two closed characteristics $(\tau,x)$ and $(\sigma,y)$ are
{\it  geometrically distinct}, if $x(\R)\not= y(\R)$. We
denote by $\J(\Sigma)$ and $\tilde\J(\Sigma)$ the set of all closed
characteristics $(\tau,x)$ on $\Sigma$ with $\tau$ being the minimal period
of $x$ and the set of all geometrically distinct ones respectively.
For $(\tau,x)\in\J(\Sigma)$, we denote by $[(\tau,x)]$ the set of all
elements in $\J(\Sigma)$ which are geometrically the same as $(\tau,x)$.
$^{\#}A$ denotes the total number of elements in a set $A$.

  To cast the given energy problem (\ref{eq1.1}) into a Hamiltonian version,
we follow \S {V.3} of  I. Ekeland's celebrated book \cite{Ek3}. Fix a
$\Sigma\in\H(2n)$ bounding a convex set $C$. Then the origin is in the
interior of $C$. Let $j_C:\R^{2n}\to [0,+\infty)$ be the gauge function of $C$
defined by
\ee
j_C(0)=0\qquad {\rm and}\qquad
j_C(x)=\inf\left\{\lambda>0\,|\,\frac{x}{\lambda}\in C\right\}
    \quad {\rm for}\quad x\ne 0.\hskip.25in \lb{eq1.2}\eee
Fix a constant $\alpha$ satisfying $1<\alpha<2$ in this paper. As usual
we define the Hamiltonian function $H_{\alpha}:\R^{2n}\to [0,+\infty)$ by
\ee
  H_{\alpha}(x) = j_C(x)^{\alpha}, \qquad \hbox{for all }x\in\R^{2n}.  \lb{eq1.3}\eee
Then $H_{\alpha}\in C^1(\R^{2n},\R)\cap C^2(\R^{2n}\bs\{0\},\R)$ is
convex and $\Sigma=H_{\alpha}^{-1}(1)$. It is well-known that
the problem (\ref{eq1.1}) is equivalent to the following given energy problem
of the Hamiltonian system
\ee
\left\{ \begin{array}{ll} \dot{x}(t) &\hskip-8pt  = JH_{\alpha}^{\prime}(x(t)),
\quad H_{\alpha}(x(t)) =1, \qquad \hbox{for all } t\in\R.\\
x(\tau) &\hskip-8pt = x(0).\end{array}\right.
\lb{eq1.4}\eee
Denote by $\J(\Sigma,\alpha)$ the set of all solutions $(\tau,x)$ of
the problem (1.4) where $\tau$ is the minimal period of
$x$, and by $\tilde{\J}(\Sigma,\alpha)$ the set of all geometrically distinct solutions of (1.4). Note that
elements in
$\J(\Sigma)$ and
$\J(\Sigma,\alpha)$ are one-to-one correspondent to each other.

Let $(\tau,x)\in \J(\Sigma,\alpha)$. The fundamental solution
$\gamma_x:[0,\tau]\to\Sp(2n)$ with $\gamma_x(0)=I$ of the linearized
Hamiltonian system
\ee \dot{y}(t) = JH_{\alpha}''(x(t))y(t), \qquad       \hbox{for all }t\in\R,  \lb{eq1.4a}\eee
is called the {\it associated symplectic path} of $(\tau,x)$.
The eigenvalues of $\gamma_x(\tau)$ are called {\it Floquet
multipliers} of $(\tau,x)$. By Proposition I.6.13 of \cite{Ek3},
the Floquet multipliers with their multiplicity and Krein signs of
$(\tau,x)\in\J(\Sigma)$ do not depend on the particular
choice of the Hamiltonian function in (1.4). For any $M\in\Sp(2n)$,
we define the {\it elliptic height} $e(M)$ of $M$ to be the total algebraic
multiplicity of all eigenvalues of $M$ on the unit circle $\U=\{z\in\C\,|\,|z|=1\}$
in the complex plane $\C$. Since $M$ is symplectic, $e(M)$ is even, and
$0\le e(M)\le 2n$. As usual \pagebreak a $(\tau,x)\in \J(\Sigma)$ is
{\it elliptic}, if
$e(\ga_x(\tau))=2n$. It is {\it nondegenerate}, if $1$ is a double Floquet
multiplier of it. It is {\it hyperbolic}, if $1$ is a double Floquet multiplier of
it, and  $e(\ga_x(\tau))=2$.  It is well-known that these concepts are independent
of the choice of $\alpha>1$.

  The study on closed characteristics in the global sense started in 1978,
when the existence of at least one closed characteristic on any
$\Sigma\in\H(2n)$ was first established by P. Rabinowitz in \cite{Ra1}
(for star-shaped hypersurfaces) and A. Weinstein in \cite{We2} independently.
In I. Ekeland and L. Lassoued's \cite{ELd}, I. Ekeland and H. Hofer's \cite{EH}
of 1987, and  A. Szulkin's \cite{Sz} of 1988, $^{\#}\tilde\J(\Sigma)\ge 2$ was
proved for any $\Sigma\in\H(2n)$ when $n\ge 2$.

  In \cite{Ek2} of  I. Ekeland in 1986 and \cite{Lo5} of Y. Long in 1998,
for any $\Sigma\in\H(2n)$ the existence of at least one nonhyperbolic
closed characteristic on $\Sigma$ was proved provided
$^{\#}\tilde\J(\Sigma)<+\infty$. In a recent paper \cite{Lo9} of Y. Long, it was
proved that for any $\Sigma\in\H(4)$, if $^{\#}\tilde\J(\Sigma)=2$, both
of the two closed characteristics must be elliptic.

  Let $\N$ denote the set of natural numbers. Our following main results in this
paper generalize the above mentioned results.

\numbereddemo{Definition}\lb{Def1.1} For $\alpha\in(1,2)$, we define a map
$\varrho_n\colon\H(2n)\to\N\cup\{ +\infty\}$
by
\ee
  \varrho_n(\Sigma)=\cases{ +\infty &if 
$^{\#}\V(\Sigma,\alpha)$\cr
&$\quad= +\infty$,\cr
  \min\{[\frac{i(x,1)+2S^+(x)-\nu(x,1)+n}{2}]\mid
      [(\tau,x)]\in\V_{\infty}(\Sigma,\alpha)\},
      &if $^{\#}\V(\Sigma,\alpha)$\cr
&$\quad< +\infty$,\cr}
   \lb{eq1.5}\eee
where $(i(x,1),\nu(x,1))$ is the Maslov-type index of $(\tau,x)$ defined in
\S 1.2, $S^+(x)$ is the splitting number of $(\tau,x)$ given in Definition~1.3, and $\V(\Sigma,\alpha)$ and
$\V_{\infty}(\Sigma,\alpha)$ are given by   Definition \ref{Def1.4} below.  \enddemo

  Note that when $^{\#}\td{\J}(\Sigma)< +\infty$, the set
$\V_{\infty}(\Sigma,\alpha)$ is nonempty and finite. By Lemma \ref{Lem5.1}
below,  $\varrho_n(\Sigma)$ does not depend on the choice of $\alpha\in(1,2)$
and is a shape invariant, i.e., is independent of dilations of $\Sg$.

\specialnumber{1.1}\proclaim{Theorem}\lb{Thm1.1} For every $\Sigma\in\H(2n)${\rm ,}  
\ee ^{\#}\tilde\J(\Sigma)\ge \varrho_n(\Sigma), \lb{eq1.6}\eee
and
\ee \varrho_n(\Sigma) \ge [\frac{n}{2}]+1, \lb{eq1.7}\eee
where $[a]=\max\{k\in\Z\,|\,k\le a\}$ for any $a\in\R${\rm .}
\et

\specialnumber{1.1}\proclaim{{C}orollary}\lb{Cor1.1} Fix $\Sigma\in\H(2n)$ and $\alpha\in (1,2)${\rm .}
Suppose every $(\tau,x)\in\J(\Sigma,\alpha)$ satisfies
\ee i(x,1)+2S^+(x)-\nu(x,1) \geq n.    \lb{eq1.8}\eee
Then 
\ee \varrho_n(\Sigma)\ge n. \lb{eq1.9}\eee
If every $(\tau,x)\in\J(\Sigma)$ is nondegenerate{\rm , (\ref{eq1.8})}
holds{\rm .} In particular{\rm ,}
\ee ^{\#}\tilde\J(\Sigma)\ge n. \lb{eq1.10}\eee
\ec

\specialnumber{1.2}\proclaim{Theorem}\lb{Thm1.2} For any $\Sigma\in\H(2n)$ satisfying
$^{\#}\tilde\J(\Sigma)< +\infty${\rm ,} there exists at least one elliptic
closed characteristic on $\Sigma${\rm .} \et

\specialnumber{1.3}\proclaim{Theorem}\lb{Thm1.3} For any $\Sigma\in\H(2n)$ satisfying
$^{\#}\tilde\J(\Sigma)< +\infty${\rm ,}  there exist at least
$\varrho_n(\Sg)-1$  {\rm (}$\;\ge [\frac{n}{2}]${\rm )} geometrically distinct closed
characteristics on $\Sigma$ possessing irrational mean indices{\rm .} \et

\specialnumber{1.4}\proclaim{Theorem}\lb{Thm1.4} Let $\Sigma\in\H(2n)$ with $n\ge 2${\rm .} Suppose
\ee  ^{\#}\tilde\J(\Sigma)\le 2\varrho_n(\Sigma)-2< +\infty. \lb{eq1.11}\eee
Then there exist at least two elliptic elements in $\tilde\J(\Sigma)${\rm .} In
particular{\rm ,} by {\rm (\ref{eq1.7})} there are at least two elliptic elements in
$\tilde\J(\Sigma)$ provided
\ee  ^{\#}\tilde\J(\Sigma)\le 2\left[\frac{n}{2}\right]. \lb{eq1.12}\eee
\et

  The study of these problems can be traced back to pioneering works
\cite{Li} of  A. Liapunov in 1892 and \cite{Hor} of V. J. Horn in 1903.
Other related significant progress  can be found in \cite{We1} of A. Weinstein,
\cite{Mo} of J. Moser, and \cite{Ba} of T. Bartsch for local results,
in \cite{Ek3} of I. Ekeland, \cite{DDE} of G. Dell'Antonio,\break B.
D'Onofrio, and I. Ekeland,  and \cite{HWZ} of H. Hofer, K. Wysocki,
and E. Zehnder  for global results,  and
in \cite{ELy} of I. Ekeland and J.-M. Lasry, \cite{AM} of A. Ambrosetti
and G. Mancini, \cite{Gi} of M. Girardi, \cite{Ho} of  H. Hofer, and
\cite{BLMR} of H. Berestycki, J.-M. Lasry, G. Mancini, B. Ruf  for results
under pinching conditions.

  A typical example of $\Sigma\in\H(2n)$ is the ellipsoid $\E_n(r)$ defined
as follows. Let $r=(r_1,\ldots,r_n)$ with $r_k>0$ for $1\leq k\leq n$. Define
\ee \E_n(r)=\left\{x=(x_1,\ldots,x_n)\in\R^{2n}\,|\,
      \frac{1}{2}\sum_{k=1}^n{{|x_k|^2}\over {r_k^2}}=1\right\}. \lb{eq1.13}\eee
If $r_j/r_k$ is irrational whenever $j\not=k$, this $\E_n(r)$ is called a
{\it weakly nonresonant ellipsoid}. In this case there are precisely
$n$ geometrically distinct closed characteristics on $\E_n(r)$,
and all of them are elliptic and nondegenerate (cf.\ \S {I.7} of \cite{Ek3}).

  It was   conjectured some time ago (cf.\ p.235 of \cite{Ek3}) that
every $\Sigma\in\H(2n)$ possesses always at least $n$ geometrically
distinct closed characteristics. By our result, we suspect that the lower
bound number $[\frac{n}{2}]+1$ found in Theorem~\ref{Thm1.1} is the best
one can hope for. We also suspect that for any $\Sigma\in\H(2n)$, if
$^{\#}\tilde\J(\Sigma)<+\infty$, every $(\tau,x)\in\J(\Sigma)$ should be
elliptic. This is true for $\H(4)$ in view of Theorem 1.6 of \cite{Lo9}
and the above Theorem \ref{Thm1.4}, and a result \cite{HWZ} by H. Hofer,
K. Wysocki, and E. Zehnder which shows that $^{\#}\tilde\J(\Sigma)<+\infty$
implies $^{\#}\tilde\J(\Sigma)=2$ when $\Sigma\in\H(4)$.

  In the rest   of this section, we introduce the quantities used in the
above theorems and describe our main ideas in their proofs. We
use some ideas from \cite{Ek3}, \cite{Lo5}, and \cite{Lo9}.

\demo{{\rm 1.2.} Maslov\/{\rm -}\/type index functions and splitting numbers} 
  As usual, the symplectic group $\Sp(2n)$ is defined by
$$ \Sp(2n) = \{M\in \GL(2n,\R)\,|\,M^TJM=J\}, $$
whose topology is the one induced from that of $\R^{4n^2}$. We are
interested in paths in $\Sp(2n)$:
$$ \P_{\tau}(2n) = \{\ga\in C([0,\tau],\Sp(2n))\,|\,\ga(0)=I\}, $$
which is equipped with the topology induced from that of  $\Sp(2n)$.
The following function is introduced in \cite{Lo6}:
$$  D_{\om}(M) = (-1)^{n-1}\ol{\om}^n\det(M-\om I), \qquad
         \hbox{for all }\om\in\U,\, M\in\Sp(2n). $$
It is proved in \cite{Lo6} that this function is real. Thus for any
$\om\in\U$ we can define
$$ \Sp(2n)_{\om}^0 = \{M\in\Sp(2n)\,|\, D_{\om}(M)=0\}.  $$
This gives a codimension-$1$ hypersurface in $\Sp(2n)$. For any
$M\in \Sp(2n)_{\om}^0$, we define a co-orientation of $\Sp(2n)_{\om}^0$
at $M$ by the positive direction $\frac{d}{dt}Me^{t\ep J}|_{t=0}$ of
the path $Me^{t\ep J}$ with $0\le t\le 1$ and $\ep>0$ sufficiently
small. We also define
\aa
\Sp(2n)_{\om}^{\ast} &=& \Sp(2n)\bs \Sp(2n)_{\om}^0,   \nn\\
\P_{\tau,\om}^{\ast}(2n) &=&
     \{\ga\in\P_{\tau}(2n)\,|\,\ga(\tau)\in\Sp(2n)_{\om}^{\ast}\}, \nn\\
\P_{\tau,\om}^0(2n) &=& \P_{\tau}(2n)\bs  \P_{\tau,\om}^{\ast}(2n).  \nn\eaa
For any two continuous arcs $\xi$ and $\eta:[0,\tau]\to\Sp(2n)$ with
$\xi(\tau)=\eta(0)$, we define as usual:
$$ \eta\ast\xi(t) = \left\{\matrix{
            \xi(2t), & \quad {\rm if}\;0\le t\le \tau/2, \cr
            \eta(2t-\tau), & \quad {\rm if}\; \tau/2\le t\le \tau. \cr}\right. $$
Given any two $2m_k\times 2m_k$ matrices of square block form
$M_k=\left(\matrix{A_k&B_k\cr
                                C_k&D_k\cr}\right)$ with $k=1, 2$,
as in \cite{DL}, the $\;\dm$-product of $M_1$ and $M_2$ is defined by
the following $2(m_1+m_2)\times 2(m_1+m_2)$ matrix $M_1\dm M_2$:
$$ M_1\dm M_2=\left(\matrix{A_1&  0&B_1&  0\cr
                              0&A_2&  0&B_2\cr
                             C_1&  0&D_1&  0\cr
                              0&C_2&  0&D_2\cr}\right). \nn$$  
Denote by $M^{\dm k}$ the $k$-fold $\dm$-product $M\dm\cdots\dm M$. Note
that the $\dm$-multiplication is associative, and the $\dm$-product of any
two symplectic matrices is symplectic. For any paths $\ga_j\in\P_{\tau}(2n_j)$
with $j=0$ and $1$, define $\ga_0\dm\ga_1(t)= \ga_0(t)\dm\ga_1(t)$ for all
$t\in [0,\tau]$.

  We define a special path $\xi_n\in\P_{\tau}(2n)$ by
\ee \xi_n(t) = \left(\matrix{2-\frac{t}{\tau} & 0 \cr
                                            0 &  (2-\frac{t}{\tau})^{-1}\cr}\right)^{\dm n},
         \quad {\rm for}\;0\le t\le \tau.  \lb{eq1.14}\eee
\enddemo

\numbereddemo{Definition}\lb{Def1.2} Let $\om\in\U$. For any $M\in \Sp(2n)$, 
\ee  \nu_{\om}(M)=\dim_{\C}\ker_{\C}(M - \om I).  \lb{eq1.15}\eee
For any $\tau>0$ and $\ga\in \P_{\tau}(2n)$,  
\ee  \nu_{\om}(\ga)= \nu_{\om}(\ga(\tau)).  \lb{eq1.15a}\eee
If $\ga\in\P_{\tau,\om}^{\ast}(2n)$, 
\ee i_{\om}(\ga) = \lrb{\Sp(2n)_{\om}^0: \ga\ast\xi_n},  \lb{eq1.16}\eee
where the right-hand side of (\ref{eq1.16}) is the usual homotopy intersection
number, and the orientation of $\ga\ast\xi_n$ is its positive time direction under
homotopy with fixed end points. If $\ga\in\P_{\tau,\om}^0(2n)$, we let $\F(\ga)$
be the set of all open neighborhoods of $\ga$ in $\P_{\tau}(2n)$, and define
\ee i_{\om}(\ga) = \sup_{U\in\F(\ga)}\inf\{i_{\om}(\beta)\,|\,
                       \beta\in U\cap\P_{\tau,\om}^{\ast}(2n)\}.
               \lb{eq1.17}\eee
Then
$$ (i_{\om}(\ga), \nu_{\om}(\ga)) \in \Z\times \{0,1,\ldots,2n\}, $$
is called the {\it index function} of $\ga$ at $\om$.
\enddemo

  Note that the right-hand side of (\ref{eq1.17}) is always finite by
Proposition 4.5 and Corollary 4.6 of \cite{Lo4}, as well as by Theorem 2.6
and Corollary 2.7 of \cite{Lo7}.

  For any symplectic path $\ga\in\P_{\tau}(2n)$ and $m\in\N$,  we
define its $m^{\rm th}$ iteration $\ga^m:[0,m\tau]\to\Sp(2n)$ by
\ee \ga^m(t) = \ga(t-j\tau)\ga(\tau)^j, \qquad
  \hbox{for all }j\tau\leq t\leq (j+1)\tau,\;j=0,1,\ldots,m-1.
     \lb{eq1.18}\eee
We still denote the extended path on $[0,+\infty)$  by $\ga$.

  Fix a $\Sigma\in\H(2n)$ and a real number $\alpha\in (1,2)$.  For any
$(\tau,x)\in\J(\Sigma,\alpha)$ and $m\in\N$, we define its $m^{\rm th}$
iteration $x^m:\R/(m\tau\Z)\to\R^{2n}$ by
\ee x^m(t) = x(t-j\tau), \qquad \hbox{for all }j\tau\leq t\leq (j+1)\tau,
       \quad j=0,1, \ldots, m-1. \lb{eq1.19}\eee
We still denote by $x$ its extension to $[0,+\infty)$.

\numbereddemo{Definition}\lb{Def1.3} For any $\ga\in\P_{\tau}(2n)$, 
\ee  (i(\ga,m), \nu(\ga,m)) = (i_1(\ga^m), \nu_1(\ga^m)), \qquad \hbox{for all }m\in\N.\hskip.5in
  \lb{eq1.19aa}\eee
The {\it mean index} $\hat{i}(\ga,m)$ per $m\tau$ for $m\in\N$ is defined by
\ee \hat{i}(\ga,m) = \lim_{k\to +\infty}\frac{i(\ga,mk)}{k}. \lb{eq1.19b}\eee
For any $M\in\Sp(2n)$ and $\om\in\U$, we define the {\it splitting numbers}
$S_M^{\pm}(\om)$ of $M$ at $\om$ by
\ee S_M^{\pm}(\om)
     = \lim_{\ep\to 0^+}i_{\om\exp(\pm\sqrt{-1}\ep)}(\ga) - i_{\om}(\ga),
   \lb{eq1.19c}\eee
for any path $\ga\in\P_{\tau}(2n)$ satisfying $\ga(\tau)=M$.

  For $\Sg\in\H(2n)$ and $\alpha\in (1,2)$, let $(\tau,x)\in \J(\Sg,\alpha)$. We define
\aa
 S^+(x) &=& S_{\ga_x(\tau)}^+(1),  \lb{eq1.19d}\\
  (i(x,m), \nu(x,m)) &=& (i(\ga_x,m), \nu(\ga_x,m)),  \lb{eq1.19e}\\
   \hat{i}(x,m) &=& \hat{i}(\ga_x,m),  \lb{eq1.19f}\eaa
for all $m\in\N$, where $\ga_x$ is the associated symplectic path of $(\tau,x)$.
\enddemo

  As proved in \cite{Lo7}, the mean index $\hat{i}(\ga,m)$ is always a finite real
number, and the splitting numbers topologically defined above are independent of
the choice of $\ga$ and a complete algebraic characterization of splitting numbers
is given by Theorem 4.11 of \cite{Lo7} (i.e., Theorem \ref{Thm6.6} below). Note
that by  (5.7) of~\cite{Lo7},  
\ee   m\hat{i}(\ga,1) = \hat{i}(\ga,m), \quad \hbox{for all }m\in\N, \ga\in\P_{\tau}(2n).
         \lb{eq1.19a}\eee

 The above Maslov-type index theory $(i_1(\ga),\nu_1(\ga))$ for $\ga\in\P_{\tau}(2n)$
was defined by C. Conley and E. Zehnder in \cite{CZ} of 1984 when $n\geq 2$ and $\ga\in\P_{\tau,1}^{\ast}(2n)$, by Y.
Long and E. Zehnder in \cite{LZ} of 1990 when $n=1$
and $\ga\in\P_{\tau,1}^{\ast}(2n)$, by Y. Long in \cite{Lo1} and C. Viterbo  \cite{Vi} in
1990 independently, when $\ga\in\P_{\tau,1}^0(2n)$ is the fundamental solution of some linear
Hamiltonian system with continuous symmetric $\tau$-periodic coefficients, and by Y. Long
in \cite{Lo4} of 1997 for any $\ga\in\P_{\tau,1}^0(2n)$. The index function
$(i_{\om}(\ga), \nu_{\om}(\ga))$ with $\om\in\U$, the Maslov-type mean index
$\hat{i}(\ga,m)$, and the splitting numbers $S_M^{\pm}(\om)$ were defined
by Y. Long in \cite{Lo7} of 1999.

\demo{{\rm 1.3.} Variational setting of the problem} 
  Fix   $\Sg\in\H(2n)$ and  $\alpha\in (1,2)$ for the rest of this section.
To solve the given fixed energy problem (1.4) as in \cite{Ek3} but
with our $J$ in (\ref{eq1.1}) instead, we consider the following fixed
period problem with $H_{\alpha}$ defined by (\ref{eq1.3}):
\ee
\left\{\matrix{\dot{z}(t) &=& JH_{\alpha}'(z(t)), & \qquad\hbox{for all }t\in\R,  \cr
                              z(1) &=& z(0). \qquad\quad &\cr}\right.
\lb{eq1.20}\eee
Define
\ee   E =\left\{u\in L^{(\alpha-1)/\alpha}(\R/\Z,\R^{2n})\,|\,\int_0^1udt=0\right\}.
\lb{eq1.21}\eee
The corresponding Clarke-Ekeland dual action functional $f:E\to\R$ is defined by
\ee f(u)=\int_0^1\lrc{\frac{1}{2}(Ju,\Pi u)+H_{\alpha}^{\ast}(-Ju)}dt,
\lb{eq1.22}\eee
where $\Pi u$ is defined by $\frac{d}{dt}\Pi u=u$ and $\int_0^1\Pi udt=0$, and the usual
dual function $H_{\alpha}^{\ast}$ of $H_{\alpha}$ is defined by
\ee  H_{\alpha}^{\ast}(x)=\sup_{y\in\R^{2n}}\left((x,y)-H_{\alpha}(y)\right).
\lb{eq1.23}\eee
Here $(\cdot,\cdot)$ denotes the standard inner product of  $\R^{2n}$. Then $f\in C^2(E,\R)$.

  Suppose $u\in E\bs\{0\}$ is a critical point of $f$. By Chapter V of \cite{Ek3}, there
exists $\xi_u\in\R^{2n}$ such that $z_u(t)=\Pi u(t)+\xi_u$ is a $1$-periodic
solution of the problem (\ref{eq1.20}). Let $h=H_{\alpha}(z_u(t))$ and $1/m$ be the
minimal period of $z_u$ for some $m\in\N$. Define
\ee
x_u(t) = h^{-1/\alpha}z_u(h^{(2-\alpha)/\alpha}t) \qquad
      {\rm and}\qquad \tau={1\over m}h^{(\alpha-2)/\alpha}.\hskip.4in
            \lb{eq1.23a}\eee
Then   $x_u(t)\in\Sigma$ for all $t\in\R$ and
$(\tau,x_u)\in\J(\Sigma,\alpha)$. Note that the  period~$1$ of $z_u$ corresponds to
the period $m\tau$ of the solution $(m\tau,x_u^m)$ of (\ref{eq1.4}) with minimal
period $\tau$.

  In \cite{Ek1} to \cite{Ek3}, I. Ekeland defined his Morse-type index theory for the functional $f$
 at its critical points $u$. The relationship between the Ekeland
index and the above Maslov-type index is given in the following lemmas.

\proclaimtitle{cf.\ \cite{Br} for the nondegenerate case, Lemma 1.3
of
\cite{Lo5} for degenerate case, and Theorem 3.2 of \cite{LZh} for a different proof}
  \specialnumber{1.1} \proclaim{Lemma}\lb{Lem1.1}  For $u$ and $z_u$ defined as  above{\rm ,}
\ee i(z_u,1)=i_1^E(u)+n \qquad {\it and}\qquad \nu(z_u,1)=\nu_1^E(u).
     \lb{eq1.23b}\eee\el

\proclaimtitle{cf.\ Lemma 1.4 of \cite{Lo5}} 
\specialnumber{1.2} \proclaim{Lemma}\lb{Lem1.2} For $z_u${\rm ,} $x_u${\rm ,}
$\tau$ and
$m$
defined as above{\rm ,}  
\ee i(x_u,m)=i(z_u,1) \qquad{\it and}\qquad \nu(x_u,m)=\nu(z_u,1).
    \lb{eq1.23c}\eee\el

  On the other hand, every solution $(\tau,x)\in\J(\Sigma,\alpha)$ gives rise to a
sequence $\{z^x_m\}_{m\in\N}$ of solutions of the given period-$1$ problem
(\ref{eq1.20}), and a sequence $\{u^x_m\}_{m\in\N}$ of critical points of $f$ defined
by  \pagebreak
\aa
  z^x_m(t) &=& (m\tau)^{-1/(2-\alpha)}x(m\tau t), \lb{eq1.23d}\\
  u^x_m(t) &=& (m\tau)^{(\alpha-1)/(2-\alpha)}\dot{x}(m\tau t).\lb{eq1.23e}
\eaa
Thus,   from the above discussion we obtain:

\proclaimtitle{cf.\ Proof of Corollary 9.4 of \cite{DL} for a
direct proof} \specialnumber{1.2}\proclaim{{C}orollary}\lb{Cor1.2} 
For any $(\tau,x)\in\J(\Sigma,\alpha)${\rm ,} 
\ee  i(x,1) \ge n.   \lb{eq1.24}\eee
\ec

  Following \S V.3 of \cite{Ek3}, we denote by "$\ind$" the Fadell-Rabinowitz
$S^1$-action cohomology index theory for $S^1$-invariant subsets of $E$
defined in \cite{Ek3} (cf.\ \cite{FR} of E. Fadell and P. Rabinowitz
for the original definition and   Appendix 2 of this paper). For
$[f]_c\equiv\{u\in E\,|\,f(u)\le c\}$, the following critical values of
$f$ are defined
\ee c_k = \inf\{c<0\,|\,\ind([f]_c)\ge k\}, \qquad \hbox{for all }k\in\N.
\lb{eq1.25}\eee
Based upon \cite{Vi1} and Proposition 2 on p. 443 of \cite{EH}, the following
important result is given in Theorem V.3.4 of \cite{Ek3}:

\specialnumber{1.5}\proclaim{Theorem}\lb{Thm1.5} For any $k\in\N${\rm ,} there exists $u\in E$ such that $f'(u)=0${\rm ,}
$f(u)=c_k${\rm ,} and
$$  i_1^E(u) \le 2k-2 \le i_1^E(u) + \nu_1^E(u) -1. $$
\et

  By the Maslov-type index theory defined above, Theorem \ref{Thm1.5}
and results of Ekeland et al.\ contained in (V.3.21), (V.3.22),
Proposition V.3.3, and Theorem V.3.4 in the Section V.3 of \cite{Ek3} can
be rephrased as follows in Theorem 1.6, which forms one of the bases of our proof.

\specialnumber{1.6}\proclaim{Theorem}\lb{Thm1.6}  
\aa
  -\infty < c_1&=& \inf_{u\in E}f(u) \le c_2 \le \cdots \le c_k\le c_{k+1}
       \le \cdots < 0,  \lb{eq1.26}\\
    c_k &\to&  0  \qquad {\rm as}\quad k\to +\infty,  \lb{eq1.26a}\\
  ^{\#}\td{\J}(\Sg) &=& +\infty \quad
{\rm if}\;\;c_k=c_{k+1}\;\;{\rm for\;some}\;\;k\in\N.
          \lb{eq1.27}\eaa
For any given $k\in\N${\rm ,} there exists $(\tau,x)\in \J(\Sg,\alpha)$ and
$m\in\N$ such that for
\aa u^x_m(t) &= &(m\tau)^{(\alpha-1)/(2-\alpha)}\dot{x}(m\tau t), \quad 0\le t\le 1, \lb{eq1.28}\\
\noalign{\noindent there hold}
  f'(u^x_m) &=& 0, \quad  f(u^x_m) = c_k,  \lb{eq1.29}\\
  i(x,m) &\le& 2k-2+n \le i(x,m) + \nu(x,m) -1.  \lb{eq1.30}\eaa
\et

  Based upon (\ref{eq1.29}) and  (\ref{eq1.30}), the following
definitions are as introduced in \cite{Lo5}.

\numbereddemo{Definition}\lb{Def1.4} For any $\Sigma\in\H(2n)$ and $\alpha\in (1,2)$,
$(\tau,x)\in\J(\Sg,\alpha)$
is $(m,k)$-{\it  variationally visible}, if there exist some $m$ and
$k\in\N$ such that (\ref{eq1.29}) and (\ref{eq1.30}) hold for $u_m^x$
defined by  (\ref{eq1.28}). We call
$(\tau,x)\in\J(\Sigma,\alpha)$ {\it infinite variationally visible}, if
there exist infinitely
many $(m,k)$ such that $(\tau,x)$ is $(m,k)$-variationally visible. We denote by
$\V(\Sg,\alpha)$ (or $\V_{\infty}(\Sg,\alpha)$) the subset of  $\td{\J}(\Sg,\alpha)$ in which a
representative
$(\tau,x)\in\J(\Sg,\alpha)$ of each $[(\tau,x)]$ is variationally visible (or infinite variationally visible). \enddemo

\demo{{\rm 1.4.} Main new ideas and sketch of  proofs} 
 We explain our ideas in the proof  of Theorem \ref{Thm1.1} first.
  As in \cite{Lo5}, we define the  $m^{\rm th}$ {\it index interval}
of  $(\tau,x)\in\J(\Sg,\alpha)$ by the closed interval
\ee \I_m(\tau,x) = [i(x,m),i(x,m)+\nu(x,m)-1].    \lb{eq1.31}\eee
We call the set
\ee \I(\tau,x)=\bigcup_{m\ge 1}\I_m(\tau,x), \lb{eq1.32}\eee
the {\it index cover set} of $(\tau,x)$. In Theorem \ref{Thm2.3} below, the
following new iteration inequality of the Maslov-type index theory is proved
for any $(\tau,x)\in\J(\Sigma,\alpha)$:
\aa\quad\qquad
i(x,m+1) - i(x,m) - \nu(x,m) &\hskip-7pt\ge\hskip-7pt&
     i(x,1) - \frac{e(\ga_x(\tau))}{2}+1    \lb{eq1.330}\\
 &\hskip-7pt\ge\hskip-7pt&  i(x,1) - n + 1,  \qquad \hbox{for all }m\in\N. \lb{eq1.33}\eaa
Here we should point out that (\ref{eq1.330}) and (\ref{eq1.33}) always hold without
the convexity condition on $\Sg$ if we delete the $1$ from the right-hand sides of
these two inequalities. To get the sharper estimate with $1$, we used the convexity
condition via the following splitting lemma.
\enddemo

\proclaimtitle{Lemma 3.2 of \cite{Lo5}}
\specialnumber{1.3} \proclaim{Lemma}\lb{Lem1.3}  Fix $\Sigma\in\H(2n)$ and
$\alpha\in (1,2)$. For any $(\tau,x)\in\J(\Sigma,\alpha)${\rm ,}
 there exist $P\in\Sp(2n)$ and $M\in\Sp(2n-2)$ such that
\ee \gamma_x(\tau) = P^{-1}(N_1(1,1)\dm M)P,   \lb{eq1.330a}\eee
where $N_1(1,1)=\left(\matrix{1&1\cr
                                                0&1\cr}\right)$.
\el

  The inequality (\ref{eq1.33}) specially implies that all index intervals of $(\tau,x)$
are mutually disjoint, and introduces a way to estimate the ellipticity of
$(\tau,x)$ in terms of its iterated Maslov-type indices for some $m\in\N$.
Note that because the system (1.4) is autonomous, we always have
\ee   \nu(x,m) \ge \nu(x,1) \ge 1.  \lb{eq1.33a}\eee
Together with Corollary \ref{Cor1.2} and (\ref{eq1.33}), we then obtain
\ee \hat{i}(x,1)\geq 2, \qquad \hbox{for all }(\tau,x) \in \J(\Sg,\alpha).  \lb{eq1.33b}\eee
\pagebreak

  Now (\ref{eq1.30}) can be restated as
\ee 2\N-2+n \subset \bigcup_{[(\tau,x)]\in\td{\J}(\Sg,\alpha)}\I(\tau,x).  \lb{eq1.34}\eee
We call integers in the sequence $2\N-2+n$ {\it effective integers}.

\vglue4pt

  Now in the following we suppose
\ee  ^{\#}\tilde\J(\Sigma)<+\infty.   \lb{eq1.36}\eee
Then by (\ref{eq1.27}) no equality in (\ref{eq1.26}) can hold, i.e., we must have
\ee -\infty < c_1 < c_2 < \cdots <c_k < c_{k+1} < \cdots < 0. \lb{eq1.38}\eee
Here we have used the multiplicity method of Fadell-Rabinowitz $S^1$-index
theory  (cf.\ \cite{FR}) via Theorem V.3.4 of \cite{Ek3}. Note that here each
$c_k$ corresponds uniquely  to an effective integer $2k-2+n$.
We observe that in this case there is a one-to-one correspondence between
the effective integers and index intervals of  all closed characteristics
$[(\tau,x)]$ in $\V_{\infty}(\Sigma,\alpha)$. In other words, under the condition
(\ref{eq1.36}), an injective map $p: \N\to \V_{\infty}(\Sigma,\alpha)\times\N$ can
be defined by  (\ref{eq1.29}), (\ref{eq1.30}), and (\ref{eq1.38}). We refer to Section \ref{s3} below for the precise
definition of the map $p$.

  From our observations on the weakly nonresonant ellipsoid as well as the
study on the case of $\H(4)$ in \cite{Lo9}, we noticed that in order to
maximize the effect of the Fadell-Rabinowitz $S^1$-index theory, instead of  the
index interval $\I_m(\tau,x)$, we should consider the largest open interval
which contains $\I_m(\tau,x)$, possesses no part of any other index interval
of $(\tau,x)$, and still can be used as the target of the map $p$. This leads
to our introduction of the index jump of $(\tau,x)$.

\numbereddemo{Definition}\lb{Def1.5} For $\Sg\in\H(2n)$ and $\alpha\in (1,2)$, we define the
{\it $m^{\rm th}$ index jump} $\G_m(\tau,x)$ of $(\tau,x)\in\J(\Sg,\alpha)$ to be
the open interval
\ee \G_m(\tau,x) = (i(x,m) + \nu(x,m) - 1,  i(x,m+2)).   \lb{eq1.40}\eee
\enddemo

   When (\ref{eq1.36}) holds,  we have $\V_{\infty}(\Sg,\alpha) \not= \emptyset$ and
we write
\ee  \V_{\infty}(\Sg,\alpha) = \{[(\tau_1,x_1)],\ldots, [(\tau_q,x_q)]\}.   \lb{eq1.42}\eee

  We will show that based on the estimates (\ref{eq1.24}) and (\ref{eq1.33b}),
there are infinitely many chances that the index jumps of all the $q$
closed characteristics contain common intervals. In fact, it will be one
of the important steps to show that there exist infinitely many
$(N,m_1,\ldots,m_q)\in\N^{q+1}$ such that
\ee \emptyset \not= [2N-\kappa_1, 2N+\kappa_2]  \subset
      \bigcap_{j=1}^q\G_{2m_j-1}(\tau_j,x_j),   \lb{eq1.44}\eee
where
\aa
\kappa_1\equiv\kappa_1(\Sigma,\alpha)&=&\min_{1\le j\le q}
  \lrp{i(x_j,1) + 2S_{\ga_x(\tau_j)}^+(1) - \nu(x_j,1)},  \lb{eq1.46}\\
\kappa_2\equiv\kappa_2(\Sigma,\alpha)&=&\min_{1\le j\le q}(i(x_j,1)-1).
    \lb{eq1.47}\eaa

  By the Fadell-Rabinowitz $S^1$-index theory and our above discussion on the map
$p$, there is a one-to-one correspondence between the effective integers contained
on the left-hand side interval of (\ref{eq1.44}) and the index jumps on the right-hand side of it. Together with
comparisons on $\kappa_1$, $\kappa_2$, and
$\varrho_n(\Sg)$, we have
\aa  q
&\ge&  ^{\#}\left( (2\N-2+n)\cap [2N-\kappa_1,2N+\kappa_2]\right) \lb{eq1.50}\\
&\ge& \varrho_n(\Sg).   \lb{eq1.52}\eaa
Then by the complete understanding on  the splitting numbers given in \cite{Lo7}
(Theorem \ref{Thm6.6} (Appendix) below), the estimate (\ref{eq1.24}), and the above
Definition \ref{Def1.1} of $\varrho_n(\Sg)$, we obtain
\ee \varrho_n(\Sg)\ge \lrb{\frac{n}{2}} + 1.  \lb{eq1.54}\eee
This yields the results of Theorem \ref{Thm1.1} when we assume (\ref{eq1.44}).

  Now the existence and size of the interval $[2N-\kappa_1,2N+\kappa_2]$ in
(\ref{eq1.44}) is very crucial for our multiplicity results.  The proof of
(\ref{eq1.44}) depends on the new abstract precise iteration formula of
the Maslov-type index theory proved in the Theorem \ref{Thm2.1} below for any
$(\tau,x)\in\J(\Sigma,\alpha)$,
\aa &&\lb{eq1.56}\\
i(x,m) &=& m\left(i(\ga,1)+S^+_M(1)- C(M)\right)
  +2\sum_{\theta\in(0,2\pi)}E\lrp{\frac{m\theta}{2\pi}}S^-_M\lrp{e^{\sqrt{-1}\theta}}   \nn\\
& &  - \lrp{S^+_M(1)+C(M)}, \quad  \hbox{for all }m\in\N, \nn\eaa
where $M=\ga_x(\tau)$, and
\aa
  E(a) &=& \min\{k\in\Z\,|\,k\ge a\},  \qquad \hbox{for all }a\in\R,  \lb{eq1.57}\\
  C(M) &=& \sum_{0<\theta<2\pi}S_M^-\lrp{e^{\sqrt{-1}\theta}}.  \lb{eq1.58}\eaa
Thus the change of $i(x,m)$ in $m$ consists of a linearly increasing term\break
$m(i(\ga,1)+S^+_M(1)- C(M))$, rotator terms $E(\frac{m\theta}{2\pi})$ with
$S_M^-(e^{\sqrt{-1}\th})>0$, and a bounded term.

  Then the control of the location and the size of the index jumps\break
$\G_{2m_j-1}(\tau_j,x_j)$ for $1\le j\le q$ depend on the control of all the
rotators in terms of the iteration time $2m_j-1$'s for $1\le j\le q$.
The corresponding rotators are divided into two sets according
to the rotation angle $\th/(2\pi)$ being rational or irrational. Now we
choose a large enough integer $m_0$ so that multiplying by $2m_0$ makes all
the rational rotation angles $2m_0\th/(2\pi)$ become integers. Then we require
that each $m_j$ has the form $d_jm_0$ for some $d_j\in\N$ to be determined later for
$1\le j\le q$. For irrational rotators, we further choose the $m_j$'s so that
\ee E\lrp{\frac{2m_j\th}{2\pi}}^{\phantom{|}}\hskip-3pt - E\lrp{\frac{(2m_j-1)\th}{2\pi}} =1 \hbox{ \rm
or }  
     E\lrp{\frac{(2m_j+1)\th}{2\pi}}- E\lrp{\frac{2m_j\th}{2\pi}}=1 \enspace \lb{eq1.60}\eee
holds. Thus the largest jumps of irrational rotators are caught. This is
realized by requiring
\ee   \lrc{\frac{m_j\th}{\pi}} \quad {\rm or}\quad 1- \lrc{\frac{m_j\th}{\pi}}
   \lb{eq1.62}\eee
 to be sufficiently small, where $\{a\}=a-[a]$ for $a\in\R$.
These requirements will imply that the index jumps
$\G_{2m_j-1}(\tau_j,x_j)$ get big enough sizes. To make them jump together
we further require that all these $m_j$'s with $1\le j\le q$ have a  common integer
factor $N\in\N$ in some sense. By choosing this $N$ carefully we fulfill the
requirements on the rational and irrational rotators, specially (\ref{eq1.62}),
simultaneously. Therefore the problem is reduced to solving the following
dynamics problem on the torus. Namely, for a given $v\in\R^k$, find infinitely
many $N\in\N$ such that the decimal part $\{Nv\}$ is as close to some vertex
$\chi$ of the cube $[0,1]^k$ as one wants; i.e., for a given small $\ep>0$,
\ee  |\{Nv\} - \chi| < \ep. \lb{eq1.64}\eee

  We observe that the closure of $\{\{Nv\}\,|\,N\in\N\}$ in the standard
torus $\T^k=\R^k/\Z^k$ forms a closed additive subgroup of $\T^k$. Thus it must
contain the identity element of $\T^k$. This proves the existence of the point
$\chi$ and infinitely many integral  $N$'s. In   Section 4 below, we prove this  in
such a way that the integers $(N,m_1,\ldots,m_q)\in\N^{q+1}$ claimed
in (\ref{eq1.44}) can be chosen simultaneously.

  Here we notice that in the above arguments the convexity of  $\Sg$ is only used
to get estimates (\ref{eq1.24}),  the splitting Lemma \ref{Lem1.3}, and the following
weaker version of (\ref{eq1.33b}),
\ee  \hat{i}(x,1) > 0, \qquad \hbox{for all }(\tau,x) \in \J(\Sg,\alpha).  \lb{eq1.66}\eee

  Our main idea in the proof of Theorem \ref{Thm1.2} is to show the existence
of one closed characteristic $[(\tau_j,x_j)]$ found by Theorem \ref{Thm1.1}
which makes   both equalities hold in (\ref{eq1.330}) and (\ref{eq1.33}) for
the chosen iteration time $m=2m_j$. Then it must be elliptic. This closed
characteristic is minimal according to the injection map $p$ in a certain sense.

  The proof of Theorem \ref{Thm1.3} depends on the understanding of  the mean
index sequence of iterations of closed characteristics.  Under the assumption
(\ref{eq1.36}),  we prove in   Lemma \ref{Lem3.1} below that according to the
ordering defined by the injection map $p$, the corresponding mean indices of
iterations of  closed characteristics strictly increase; i.e.,
\ee  0 <  \hat{i}(x_{j(s)},2m_{j(s)}) < \hat{i}(x_{j(t)},2m_{j(t)}),
     \lb{eq1.68}\eee
for $p(k) = ([(\tau_{j(k)},x_{j(k)})], 2m_{j(k)})$, with $k=s$ or $t$, and
$1\le s<t\le \varrho_n(\Sg)$.  Then we prove that if both  mean
indices $\hat{i}(x_{j(s)},1)$ and $\hat{i}(x_{j(t)},1)$ are rational, by our choice
of the iteration time $m_j$'s, the two iterated mean indices in (\ref{eq1.68}) must
be equal to each other. This yields a contradiction.

  To prove Theorem \ref{Thm1.4}, we further observe that the elliptic solution
found in   Theorem \ref{Thm1.2} corresponds to the vertex $\chi$ of the cube
$[0,1]^k$ in (\ref{eq1.64}) via the injection map $p$. By Theorem \ref{Thm1.3},
when $n\ge 2$ there exist at least two such vertices which make (\ref{eq1.64}) hold.
Then we prove that they produce two different elliptic closed orbits.

  This paper is organized as follows. In Section~\ref{s2}, we derive the abstract
precise iteration formula (\ref{eq1.56}) and the iteration inequality
(\ref{eq1.33}). In Section~\ref{s3}, we give the precise definition of the injection
map $p$ from effective numbers to iterations of closed characteristics.
In Section~\ref{s4}, we prove the common index jump theorem based on properties
of the torus group and iteration properties of the Maslov-type index theory
established in \cite{Lo7} and \cite{Lo9}. In Section~\ref{s5}, we give the proofs of
  Theorems \ref{Thm1.1} to \ref{Thm1.4}. For the reader's convenience, we give
a brief review on the Maslov-type index and its iteration theory in Section~\ref{s6}.
In Section~\ref{s7}, an appendix on the Fadell-Rabinowitz cohomology index given
by John Mather is included.

\section{Iteration formula and inequalities\\ of the Maslov-type index theory}\lb{s2}

  We refer readers to Sections~\ref{s1} and \ref{s6} for a brief review on the
Maslov-type index theory and its iteration theory. In the following,  we shall
establish an abstract precise iteration formula and new iteration inequalities
for the Maslov-type index theory using notation  in those two sections.

\demo{{\rm 2.1.} Abstract precise iteration formulae of Maslov-type indices} 
 For any $M\in\Sp(2n)$, by Lemma 4.6 of \cite{Lo7} on the splitting
numbers defined by (\ref{eq1.19c}), $S_M^-(\om)=0$ if $\om\not\in\sg(M)$.
Thus   $C(M)=\sum_{0<\theta<2\pi}S_M^-(e^{\sqrt{-1}\theta})$
defined in (\ref{eq1.58}) is a finite sum. For any $x\in\R$, in addition to the
function $E(x)$ given by (\ref{eq1.57}), we further define functions
$\,[\,\cdot\,], \;\;\phi(\cdot) : \R\to \Z$ and $ \{\,\cdot\,\}:\R\to (0,1)$ by
\y
  [x] &=& \max\{k\in\Z\mid k\le x\},\\
  \phi(x) &=& E(x) - [x], \\
  \{x\}& =&x-[x].     \ey
Note particularly that $\phi(x)=0$ if $x\in\Z$, and $\phi(x)=1$ if $x\not\in \Z$.

  Motivated by the precise iteration formula Theorem 1.3 of  \cite{Lo9}
(Theorem \ref{Thm6.7} below),  we prove the following abstract precise
iteration formula:
\enddemo

\specialnumber{2.1}\proclaim{Theorem}\lb{Thm2.1} For $n\in\N$, $\tau>0$ and any path $\ga\in\P_{\tau}(2n)${\rm ,}
set $M=\ga(\tau)$. Extend $\ga$ to the whole $[0,+\infty)$ by  {\rm (\ref{eq1.18}).}
Then for any  $m\in\N${\rm ,}
\aa
\qquad i(\ga,m) &=& m(i(\ga,1)+S^+_M(1)-C(M)) 
\lb{eq2.12}\\
& & +\ 2\sum_{\theta\in(0,2\pi)}E\left(\frac{m\theta}{2\pi}\right)S^-_M\left(e^{\sqrt{-1}\theta}\right)    -
\lrp{S^+_M(1)+C(M)}.\nn\eaa
\et

  \demo{Proof} Note that by Section 4 of \cite{Lo7}, for a fixed
path $\ga$ the index $i_{\omega}(\ga)$ is a step function in $\omega\in\U$
with possible jumps only at eigenvalues of $M=\ga(\tau)$ on $\U$. The
splitting numbers $S_M^{\pm}(\omega_0)$ measure the jumps between
$i_{\omega_0}(\ga)$ and nearby $i_{\omega}(\ga)$ from two sides of
$\omega_0$ in $\U$. There,    $S_M^{\pm}(\omega)=0$ if
$\omega\not\in\sigma(M)$. Therefore for any
$\omega_0=e^{\sqrt{-1}\theta_0}\in\U$ with $0\le \theta_0<2\pi$, we denote
by $\omega_j$, with $1\le j\le p_0$, the eigenvalues of $M$ on $\U$ which
are distributed counterclockwise  from $1$ to $\om_0$ and located strictly
between $1$ and $\omega_0$. Then we have
\aa
i_{\omega_0}(\ga)
&=& i(\ga,1) + S_M^+(1) + \sum_{j=1}^{p_0}\lrp{-S_M^-(\om_j) + S_M^+(\om_j)}
                     - S_M^-(\om_0) \nn\\
&=& i(\ga,1) + \sum_{0\le\theta<\theta_0}S^+_M\lrp{e^{\sqrt{-1}\theta}}
     - \sum_{0<\theta\le \theta_0}S^-_M\lrp{e^{\sqrt{-1}\theta}}. \nn
\eaa
Thus by the Bott-type formula in Appendix 1, Theorem \ref{Thm6.3}, for
any $m\in\N$,
\aa
&&\lb{eq2.13}\\
i(\ga,m)
&\hskip-8pt=\hskip-8pt& \sum_{\omega^m=1}i_{\omega}(\ga) \nn\\
&\hskip-8pt=\hskip-8pt& i(\ga,1) + \sum_{k=1}^{m-1}\left(i(\ga,1)
    +\hskip-4pt \sum_{0\le\theta<\frac{2k\pi}{m}}S^+_M\lrp{e^{\sqrt{-1}\theta}}
       -\hskip-4pt \sum_{0<\theta\le\frac{2k\pi}{m}}S^-_M\lrp{e^{\sqrt{-1}\theta}}\hskip-4pt \right)  
\nn\\ &\hskip-8pt =\hskip-8pt& i(\ga,1) + (m-1)\lrp{i(\ga,1) + S^+_M(1)} \nn\\
&\hskip-8pt \hskip-8pt &  + \sum_{\theta\in(0,2\pi)}\left(
       \sum_{\frac{m\theta}{2\pi}<k\le m-1}S^+_M\lrp{e^{\sqrt{-1}\theta}}
      -\sum_{\frac{m\theta}{2\pi}\le k\le m-1}S^-_M\lrp{e^{\sqrt{-1}\theta}}\right)
                      \nn\\
&\hskip-8pt=\hskip-8pt& i(\ga,1) + (m-1)\lrp{i(\ga,1) + S^+_M(1)}    \lb{eq2.14}\\
&\hskip-8pt \hskip-8pt &  +
\sum_{\theta\in(0,2\pi)}\left(\lrp{m-1-\lrb{\frac{m\theta}{2\pi}}}S^+_M\lrp{e^{\sqrt{-1}\theta}}\right.
\nn\\
&\hskip-16pt & 
      \left.  \hskip.75in  -\lrb{\frac{m(2\pi-\theta)}{2\pi}}S^-_M\lrp{e^{\sqrt{-1}\theta}}\right). \nn\eaa
Here to get  (\ref{eq2.14}) from (\ref{eq2.13}), we have used the fact that  $[m-x]=m-[x]$ for
$x\in\Z$ and $[m-x]=m-[x]-1$ for $x\not\in\Z$ to count the number of $k$'s in the
given intervals.

  By  Lemma 4.6 of \cite{Lo7}, 
\ee  S_M^+(\om)=S_M^-(\ol{\om}), \qquad \hbox{for all }\om\in\U.  \lb{eq2.14a}\eee
Thus by regrouping terms in (\ref{eq2.14}) according to $S_M^-(\omega)$ with
$\omega\in\U$, we get
\aa
&& \lb{eq2.15}\\
i(\ga,m)
&=& i(\ga,1) + (m-1)\lrp{i(\ga,1) + S^+_M(1)} + \lrp{m-1-2\lrb{\frac{m}{2}}}S_M^-(-1) \nn\\
& &  +
\sum_{\theta\in(0,\pi)}\left(\lrp{m-1-\lrb{\frac{m(2\pi-\theta)}{2\pi}}}S^-_M\lrp{e^{\sqrt{-1}\theta}}
                      \right.\nn\\
&&\left.\hskip.75in + \lrp{m-1-\lrb{\frac{m\theta}{2\pi}}}S^-_M\lrp{e^{\sqrt{-1}(2\pi-\theta)}} \right)
\nn\\ & &  - \sum_{\theta\in(0,\pi)}\left(\lrb{\frac{m(2\pi-\theta)}{2\pi}}S^-_M\lrp{e^{\sqrt{-1}\theta}}
                       +\lrb{\frac{m\theta}{2\pi}}S^-_M\lrp{e^{\sqrt{-1}(2\pi-\theta)}}\right) \nn\\
&=& i(\ga,1) + (m-1)\lrp{i(\ga,1) + S^+_M(1)} + \lrp{m-2\lrb{\frac{m}{2}}-1}S^-_M(-1) \nn\\
& &  +\sum_{\theta\in(0,\pi)}\left(\lrp{m-1-2\lrb{\frac{m\theta}{2\pi}}}
                 S^-_M\lrp{e^{\sqrt{-1}(2\pi-\theta)}}\right. \nn\\
& &\;\;\qquad\quad\left. +\lrp{m-1-2\lrb{\frac{m(2\pi-\theta)}{2\pi}}}S^-_M\lrp{e^{\sqrt{-
      1}\theta}}\right).         \nn\eaa
Since $E(x)+[y]=x+y$ if $x,y\in\R$ and $x+y\in\Z$, from (\ref{eq2.15}) we obtain
\aa
\qquad i(\ga,m)\hskip-4pt
&=&\hskip-4pt i(\ga,1) + (m-1)\lrp{i(\ga,1) + S^+_M(1)}\lb{eq2.16}\\
&&\hskip-4pt + \lrp{2E\lrp{\frac{m}{2}}-(m+1)}S^-_M(-1) \nn\\
& &\hskip-4pt +\sum_{\theta\in(0,\pi)}\left(\lrp{2E\lrp{\frac{m(2\pi-\theta)}{2\pi}}-(m+1)}
              S^-_M\lrp{e^{\sqrt{-1}(2\pi-\theta)}}\right.\nn\\
& &\hskip-4pt \qquad \quad\left.+\lrp{2E\lrp{\frac{m\theta}{2\pi}}-(m+1)}S^-_M\lrp{e^{\sqrt{-
         1}\theta}}\right)\nn\\
&=&\hskip-4pt i(\ga,1) + (m-1)\lrp{i(\ga,1) + S^+_M(1)} - (m+1)C(M)  \nn\\
& &\hskip-4pt + \sum_{\theta\in(0,2\pi)}2E\lrp{\frac{m\theta}{2\pi}}S^-_M\lrp{e^{\sqrt{-1}\theta}} .
    \nn\eaa
Here to get (\ref{eq2.16}), we have used the definition (\ref{eq1.58}) of $C(M)$.
This yields (\ref{eq2.12}). \enddemo

\demo{{R}emark {\rm 2.1}} By direct verification, Theorem \ref{Thm2.1} coincides with Theorem
\ref{Thm6.7} below. \er

\proclaimtitle{cf.\ Theorem 1.5 of \cite{Lo7}, Corollary 5.1 of
\cite{Lo9}}
\specialnumber{2.1}\proclaim{{C}orollary}\lb{Cor2.1}  For any $\tau>0$ and $\ga\in\P_{\tau}(2n)${\rm ,} 
set $M=\ga(\tau)$. Now{\rm ,}
\ee
\hat{i}(\ga,1) \equiv \lim_{k\to +\infty}\frac{i(\ga,k)}{k}
   =  i(\ga,1) + S^+_M(1) - C(M)
       +\sum_{\theta\in(0,2\pi)}\frac{\theta}{\pi}S^-_M\lrp{e^{\sqrt{-1}\theta}}.\;
        \lb{2.16a}\eee
\ec

\demo{{\rm 2.2.} New iteration inequalities of Maslov-type indices} 
  Now we give new increasing estimates of the iterated Maslov-type indices
by the following Theorems 2.2, 2.3 and 2.4. The proof of our main results
in later sections of this paper will depend only on   Theorems 2.2 and
2.3 here. The proof of Theorem 2.2 relies on the precise iteration formula
in Theorem \ref{Thm6.7} below. This method actually gives us a way to detect
and to prove or disprove whether a proposed iteration inequality of the
Maslov-type index theory is correct whenever it can be reduced to an
inequality of symplectic additive index terms. Theorem 2.4 is a
generalization of Theorem 2.2. Our proof of Theorem 2.4 is rather
different from and shorter than that of Theorem 2.2. This proof depends
on special properties of the proposed inequality, the above
Theorem \ref{Thm2.1}, and   complete understanding of splitting
numbers given by Theorem \ref{Thm6.5} below.

\specialnumber{2.2}\proclaim{Theorem}\lb{Thm2.2} For $n\in\N${\rm ,} $\tau>0${\rm ,} and
$\ga\in\P_{\tau}(2n)${\rm ,} set  $M=\ga(\tau)${\rm .}  We extend $\ga$ to $[0,+\infty)$ by {\rm (\ref{eq1.18}).}
Then for any $m\in\N${\rm ,} there holds
\aa
\nu(\ga,m) - \frac{e(M)}{2} &\le& i(\ga,m+1) - i(\ga,m) - i(\ga,1) \lb{eq2.20}\\
    &\le& \nu(\ga,1) - \nu(\ga,m+1) + \frac{e(M)}{2},  \lb{eq2.21}\eaa
where  $e(M)$ is the elliptic height defined in Section  {\rm \ref{s1}.}
\et

\demo{Proof} Without loss of generality, we may set $\tau=1$.
  By Theorem \ref{Thm6.5}, there exist $f\in C([0,1],\Omega^0(M))$ and basic
normal forms $M_1,\ldots,M_p\in\Sp(2)$ and $N_1,\ldots,N_q\in\Sp(4)$ such that
\ee   f(0)=M,  \quad
  f(1)= M_1\dm \cdots \dm M_p\dm N_1\dm \cdots\dm  N_q.  \lb{eq2.23}\eee
Since $\Sp(2n)$ is path connected, we can pick up paths $\xi_i\in \P_1(2)$ and
$\eta_j\in\P_1(4)$ such that $\xi_i(1)=M_i$ and $\eta_j(1)=N_j$ for
$1\le i\le p$ and $1\le j\le q$.

  For $k\in\Z$, we define a new path $\beta\in\P_1(2n)$ by
\ee  \beta = \left\{\matrix{
    (\xi_1\ast\phi_{2k\pi,1})\dm \xi_2\dm\cdots\xi_p\dm\eta_1\dm\cdots\dm\eta_q,
             & \;\;\mbox{if}\;\;p\geq 1, \cr
    (\eta_1\ast(\phi_{2k\pi,1}\dm I_2))\dm \eta_2\dm\cdots\dm\eta_q,
             &\;\;\mbox{if}\;\;p=0, \cr}\right.  \hskip.5in \lb{eq2.25}\eee
where $\phi_{\th,1}$ is as defined in Section~\ref{s6}. By this construction and the definition of $f(1)$,  
\ee \nu(\beta,1) = \nu(\ga,1) \quad \mbox{and}\quad e(\beta(1))=e(f(1))\le e(M).
      \lb{eq2.26}\eee
Since the curve $f$ is in $\Omega^0(M)$, we can fix a $k\in\Z$ in (\ref{eq2.25}) such that
\ee  i(\beta,1) = i(\ga,1). \lb{eq2.27}\eee
Thus by Theorem \ref{Thm6.2},
\ee  \beta\sim_1\ga \quad \mbox{on}\;\;[0,1]\;\;\mbox{along}\quad f.  \lb{eq2.28}\eee
Extending this homotopy to $[0,1]\times [0,m]$ for any $m\in\N$, we have
\ee \beta\sim_1\ga \quad \mbox{on}\;\;[0,m]\;\;\mbox{along}\quad f^m,
\lb{eq2.29}\eee
where $f^m(s)=f(s)^m$ for any $s\in [0,1]$. Then by Theorem \ref{Thm6.1},
we obtain
\ee  i(\beta,m) = i(\ga,m), \quad  \nu(\beta,m) = \nu(\ga,m), \quad \hbox{for all }m\in\N. 
\hskip.5in\lb{eq2.30}\eee
So by (\ref{eq2.26}) and (\ref{eq2.30}), it suffices to prove (\ref{eq2.20}) and (\ref{eq2.21}) for the path $\beta$.

  Note that all terms in the two inequalities (\ref{eq2.20}) and (\ref{eq2.21}),
except the elliptic height $e(\ga(1))$, are symplectically additive and homotopy invariant. When we consider the path
$\beta$, the elliptic height is also additive
in terms of the decomposition of $\beta(1)=f(1)$ in (\ref{eq2.23}). Thus by
(\ref{eq2.23}), (\ref{eq2.25}), and Theorem \ref{Thm6.1}, the proofs of both
(\ref{eq2.20}) and (\ref{eq2.21}) for the path $\beta$ are reduced to
those for each component path of  $\beta$ in (\ref{eq2.25}), i.e., any  path in
$\P_1(2)$ or $\P_1(4)$ ending at a basic normal form of symplectic matrices defined
in Section~\ref{s6}.

  Next we continue our
 proof on paths ending at these normal forms in eight cases according to all the different patterns of iteration
formulae.

\demo{Case 1}  $\ga\in\P_{\tau}(2)$ {\it and} $M\equiv \ga(1)=N_1(1,b)$ {\it with} $b=1$ {\it or} $0$. 
  In this case by Theorem \ref{Thm6.7} (or Theorem 3.4 of \cite{Lo9}), we have
\aa
 \qquad\qquad  e(M) &=& 2, \quad \nu(\ga,1) = 2-b, \lb{eq2.41}\\
    i(\ga,m) &=& m(i(\ga,1)+1)-1, \quad \nu(\ga,m)=\nu(\ga,1), \quad \hbox{for all }m\in\N. \lb{eq2.42}\eaa
Thus,
\aa
\nu(\ga,m) - \frac{e(M)}{2} &=& 1-b,  \nn\\
i(\ga,m+1) - i(\ga,m) - i(\ga,1) &=& 1, \nn\\
\nu(\ga,1) - \nu(\ga,m+1) + \frac{e(M)}{2} &=& 1. \nn\eaa
Then both (\ref{eq2.20}) and (\ref{eq2.21}) hold.

  Here we notice that when $M=N_1(1,1)$, the left-hand side of the inequality (\ref{eq2.20})
can be improved by $1$. We shall see this point in 
 Theorem \ref{Thm2.3} below.
\enddemo

\demo{Case 2} \hskip-10pt
{\it $\ga\in\P_{\tau}(2)$ and $M\equiv \ga(1)=N_1(1,-1)$.} 
 \hskip-2pt In this case by Theorem~\ref{Thm6.7} (or Theorem 3.4 of \cite{Lo9}), we have
\aa
    e(M) &=& 2, \quad \nu(\ga,1) = 1, \lb{eq2.43}\\
    i(\ga,m) &=& m i(\ga,1), \quad \nu(\ga,m)=\nu(\ga,1), \quad \hbox{for all }m\in\N. \lb{eq2.44}\eaa
Thus,
\aa
\nu(\ga,m) - \frac{e(M)}{2} &=& 0,  \nn\\
i(\ga,m+1) - i(\ga,m) - i(\ga,1) &=& 0, \nn\\
\nu(\ga,1) - \nu(\ga,m+1) + \frac{e(M)}{2} &=& 1. \nn\eaa
Then both (\ref{eq2.20}) and (\ref{eq2.21}) hold.
\enddemo

\demo{Case 3} {\it $\ga\in\P_{\tau}(2)$ and $M\equiv \ga(1)=N_1(-1,-b)$ with $b=1$ or $0$.}
  In this case by Theorem \ref{Thm6.7} (or Theorem 3.5 of \cite{Lo9}), we have
\aa
  e(M) &=& 2, \quad \nu(\ga,1) = 0, \lb{eq2.45}\\
i(\ga,m) &=& m i(\ga,1)-\frac{1+(-1)^m}{2}, \lb{eq2.46}\\
\nu(\ga,m) &=&\frac{1+(-1)^m}{1+b},
   \quad         \hbox{ for all }m\in\N. \nn\eaa
Thus,
\aa
\nu(\ga,m) - \frac{e(M)}{2} &=& \frac{(-1)^m-b}{1+b},  \nn\\
i(\ga,m+1) - i(\ga,m) - i(\ga,1) &=& (-1)^m, \nn\\
\nu(\ga,1) - \nu(\ga,m+1) + \frac{e(M)}{2} &=& \frac{(-1)^m+b}{1+b}. \nn\eaa
Then both (\ref{eq2.20}) and (\ref{eq2.21}) hold.
\enddemo

\demo{Case 4} {\it $\ga\in\P_{\tau}(2)$ and $M\equiv \ga(1)=N_1(-1,1)$.}
In this case by Theorem~\ref{Thm6.7} (or Theorem 3.5 of \cite{Lo9}), we have
\aa
    e(M) &=& 2, \quad \nu(\ga,1) = 0, \lb{eq2.49}\\
    i(\ga,m) &=& m i(\ga,1), \quad \nu(\ga,m)= \frac{1+(-1)^m}{2}, \quad
            \hbox{for all }m\in\N. \lb{eq2.50}\eaa
Thus,
\aa
\nu(\ga,m) - \frac{e(M)}{2} &=& \frac{(-1)^m-1}{2} \le 0,  \nn\\
i(\ga,m+1) - i(\ga,m) - i(\ga,1) &=& 0, \nn\\
\nu(\ga,1) - \nu(\ga,m+1) + \frac{e(M)}{2} &=& \frac{1-(-1)^m}{2} \ge 0. \nn\eaa
Then both (\ref{eq2.20}) and (\ref{eq2.21}) hold.
\enddemo

\demo{Case 5} {\it $\ga\in\P_{\tau}(2)$ and $M\equiv \ga(1)=D(2)$ or $D(-2)$.}
In this case by Theorem \ref{Thm6.7} (or Theorem 3.6 of \cite{Lo9}), we have
\aa
    e(M) &=& 0, \quad \nu(\ga,1) = 0, \lb{eq2.51}\\
    i(\ga,m) &=& m i(\ga,1), \quad \nu(\ga,m)= 0, \quad
            \hbox{for all }m\in\N. \lb{eq2.52}\eaa
Thus,
\aa
\nu(\ga,m) - \frac{e(M)}{2} &=& 0,  \nn\\
i(\ga,m+1) - i(\ga,m) - i(\ga,1) &=& 0, \nn\\
\nu(\ga,1) - \nu(\ga,m+1) + \frac{e(M)}{2} &=& 0. \nn\eaa
Then both (\ref{eq2.20}) and (\ref{eq2.21}) hold.
\enddemo

\demo{Case 6} {\it $\ga\in\P_{\tau}(2)$ and $M\equiv\ga(1)=R(\th)$ for some
$\th\in (0,\pi)\cup (\pi,2\pi)$.}
  In this case by Theorem \ref{Thm6.7} (or Theorem 3.7 of \cite{Lo9}), we have
\aa
   e(M) &=& 2, \quad \nu(\ga,1) = 0, \lb{eq2.53}\\
    i(\ga,m) &=& m \lrp{i(\ga,1)-1}+2E\lrp{\frac{m\th}{2\pi}}-1,\lb{eq2.54} \\
\nu(\ga,m)&=&2-2\phi\lrp{\frac{m\th}{2\pi}}, \;
            \hbox{for all }m\in\N.\nn \eaa
Thus,
\aa
\nu(\ga,m) - \frac{e(M)}{2} &\hskip-7pt=\hskip-7pt& 1-2\phi\lrp{\frac{m\th}{2\pi}},  \lb{eq2.55}\\
&& \lb{eq2.56}\\
i(\ga,m+1) - i(\ga,m) - i(\ga,1) &\hskip-7pt=\hskip-7pt& 2\left(E\lrp{\frac{(m+1)\th}{2\pi}}-E\lrp{\frac{m\th}{2\pi}}\right)
-1,
          \nn\\
\quad\qquad\nu(\ga,1) - \nu(\ga,m+1) + \frac{e(M)}{2} &\hskip-7pt=\hskip-7pt& 2\phi\lrp{\frac{(m+1)\th}{2\pi}}-1.
\lb{eq2.57}\eaa Note that 
\aa
&&\hskip-40pt 2\left(E\lrp{\frac{(m+1)\th}{2\pi}}
-E\lrp{\frac{m\th}{2\pi}}\right) - 1 \lb{eq2.58a}\\
  &=& 1-2\phi\lrp{\frac{m\th}{2\pi}}
           + 2\left(E\lrp{\frac{(m+1)\th}{2\pi}}-\lrb{\frac{m\th}{2\pi}}-1\right) \nn\\
  &=& 2\phi\lrp{\frac{(m+1)\th}{2\pi}}-1 - 2\left(E\lrp{\frac{m\th}{2\pi}}
 - \lrb{\frac{(m+1)\th}{2\pi}}\right),
             \lb{eq2.59a}\eaa
and that 
\aa
  E\lrp{\frac{(m+1)\th}{2\pi}} - \lrb{\frac{m\th}{2\pi}}-1  &\ge&  0, \lb{eq2.60a}\\
  E\lrp{\frac{m\th}{2\pi}} - \lrb{\frac{(m+1)\th}{2\pi}} &\ge&  0. \lb{eq2.60b}\eaa
With  (\ref{eq2.55})--(\ref{eq2.59a}), we obtain (\ref{eq2.20}) and
(\ref{eq2.21}).
\enddemo

\demo{Case 7} {\it $\ga\in\P_{\tau}(4)$ and $M\equiv\ga(1)=N_2(\om,b)$ for
some $\om\in\U\bs\R$ is nontrivial.}
  In this case by Theorem \ref{Thm6.7} (or Theorem 4.5 of \cite{Lo9}), we have
\aa
 e(M) &=& 4, \quad \nu(\ga,1) = 0, \lb{eq2.58}\\
 i(\ga,m) &=& m i(\ga,1) +2\phi\lrp{\frac{m\th}{2\pi}}  -2, \lb{eq2.59} \\
 \nu(\ga,m)
& =& 2-2\phi\lrp{\frac{m\th}{2\pi}},
              \quad    \hbox{for all }m\in\N.\nn\eaa
Thus,
\aa
\nu(\ga,m) - \frac{e(M)}{2} &=& -2\phi\lrp{\frac{m\th}{2\pi}},  \nn\\
i(\ga,m+1) - i(\ga,m) - i(\ga,1) &=& 2\phi\lrp{\frac{(m+1)\th}{2\pi}}-
  2\phi\lrp{\frac{m\th}{2\pi}},  \nn\\
\nu(\ga,1) - \nu(\ga,m+1) + \frac{e(M)}{2} &=& 2\phi\lrp{\frac{(m+1)\th}{2\pi}}.
 \nn\eaa
Then both (\ref{eq2.20}) and (\ref{eq2.21}) hold.
\enddemo

\demo{Case 8} {\it $\ga\in\P_{\tau}(4)$ and $M\equiv\ga(1)=N_2(\om,b)$ for
some $\om\in\U\bs\R$ is trivial.}
 In this case by Theorem \ref{Thm6.7} (or Theorem 4.6 of \cite{Lo9}), we have
\aa
 e(M) &=& 4, \quad \nu(\ga,1) = 0, \lb{eq2.60}\\
 \quad i(\ga,m) &=& m i(\ga,1), \quad \nu(\ga,m)= 2-2\phi\lrp{\frac{m\th}{2\pi}}, \quad
            \hbox{for all }m\in\N. \lb{eq2.61}\eaa
Thus,
\aa
\nu(\ga,m) - \frac{e(M)}{2} &=& -2\phi\lrp{\frac{m\th}{2\pi}},  \nn\\
i(\ga,m+1) - i(\ga,m) - i(\ga,1) &=& 0, \nn\\
\nu(\ga,1) - \nu(\ga,m+1) + \frac{e(M)}{2} &=& 2\phi\lrp{\frac{(m+1)\th}{2\pi}}. \nn\eaa
Then both (\ref{eq2.20}) and (\ref{eq2.21}) hold.
 The proof is complete.\hfill\qed \enddemo

  Suggested by the study of  Case 1 with $b=1$ of the proof of Theorem
\ref{Thm2.1}, we have the following result. Note that by Lemma \ref{Lem1.3}
such a consideration is useful for the study of closed characteristics on
convex hypersurfaces in $\R^{2n}$.

\specialnumber{2.3}\proclaim{Theorem}\lb{Thm2.3} For $n\in\N$, $\tau>0${\rm ,} and $\ga\in\P_{\tau}(2n)${\rm ,}
set  $M=\ga(\tau)${\rm .} Extend $\ga$ to $[0,+\infty)$ by {\rm (\ref{eq1.18}).}
Suppose that there exist $P\in\Sp(2n)$ and $Q\in\Sp(2n-2)$ such that
\ee  M=P^{-1}(N_1(1,1)\dm Q)P.  \lb{eq2.70}\eee
Then for any $m\in\N${\rm ,} 
\aa
\nu(\ga,m) - \frac{e(M)}{2}+1 &\le& i(\ga,m+1) - i(\ga,m) - i(\ga,1) \lb{eq2.71}\\   &\le& \nu(\ga,1) - \nu(\ga,m+1) +
\frac{e(M)}{2}.  \lb{eq2.72}\eaa
\et

\demo{Proof}  Without loss of generality, we may  set $\tau=1$.
  Using the path $\phi_{\th,1}$ given in Section~\ref{s6}, for $k\in\Z$ we define
\ee  \zeta(t) = N_1(1,t), \quad \xi(t)=\zeta\ast\phi_{2k\pi,1}(t), \quad
               \hbox{for all }t\in [0,1].\hskip.4in \lb{eq2.73}\eee
Then by Theorems \ref{Thm6.1} and \ref{Thm6.7} below,
\aa
  i(\xi,1) &=& 2k-1,          \lb{eq2.74}\\
  i(\xi,m) &=& 2mk - 1, \quad \nu(\xi,m)=1, \quad \hbox{for all }m\in\N,   \lb{eq2.75}\eaa
Since $\Sp(2n-2)$ is path connected, there is a path $\eta\in\P_1(2n-2)$ such that
$\eta(1)=Q$. By (\ref{eq2.70}), we can fix a $k\in\Z$ so that
\ee  2k-1+i(\eta,1) = i(\xi,1) + i(\eta,1) = i(\ga,1).  \lb{eq2.77}\eee
Pick   a path $g\in\P_1(2n)$ such that $g(1)=P$. Define
$$ h(s) = g(s)^{-1}Mg(s), \qquad \hbox{for all }t\in [0,1]. $$
Then $h$ connects $M=\ga(1)$ to $N_1(1,1)\dm Q=(\xi\dm\eta)(1)$ within $\Om^0(M)$.
Thus by Theorem \ref{Thm6.2} below,  
$$  \ga\sim_1 (\xi\dm \eta)\quad {\rm on} \quad [0,1]\quad {\rm along }\quad h.  $$
By extending this homotopy  map to $[0,1]\times [0,m]$ via iteration, we obtain
$$  \ga\sim_1 (\xi\dm\eta) \quad {\rm on}\quad [0,m]\quad {\rm along}\quad h^m, $$
where $h^m(s)=h(s)^m$. Therefore by  Theorem \ref{Thm6.1} below
and (\ref{eq2.75}), for any $m\in\N$ we obtain
\aa
  e(M) &=& e(N_1(1,1)) + e(Q) = 2 + e(Q), \lb{eq2.78}\\
  i(\ga,m) &=& i(\xi,m) + i(\eta,m) = 2mk-1 + i(\eta,m),   \lb{eq2.79}\\
  \nu(\ga,m) &=& \nu(\xi,m) + \nu(\eta,m) = 1 + \nu(\eta,m). \lb{eq2.80}\eaa
Thus plugging (\ref{eq2.78})--(\ref{eq2.80}) into (\ref{eq2.71}) and (\ref{eq2.72}),
they  become
\aa
\nu(\eta,m)  - \frac{e(Q)}{2} &\le& i(\eta,m+1)  - i(\eta,m)  - i(\eta,1) \nn\\
    &\le& \nu(\eta,1) - \nu(\eta,m+1) + \frac{e(Q)}{2},  \nn\eaa
which follows from Theorem \ref{Thm2.2} for the path $\eta$. \enddemo

  We start our generalization of Theorem \ref{Thm2.2} from the following
two lemmas.

\specialnumber{2.1} \proclaim{Lemma}\lb{Lem2.1} For any $\om\in\U$ and $M\in\Sp(2n)${\rm ,}   denote by
$(p_{\om}(M),q_{\om}(M))$ the Krein type of $\om\in\sg(M)${\rm .} Then
\aa
   0 &\le& \nu_{\om}(M) - S_M^-(\om)\le p_{\om}(M),  \lb{eq2.81}\\
   0 &\le& \nu_{\om}(M) - S_M^+(\om)\le q_{\om}(M).  \lb{eq2.82}\eaa
\el

\demo{Proof} The left inequalities are proved by Corollary 4.13 of \cite{Lo7}.
Denote by $N$ the matrix on the right-hand side of (\ref{eq6.10})
in Appendix 1, Theorem 6.5. Then
$$ p_{\om}(N)\le p_{\om}(M), \quad q_{\om}(N)\le q_{\om}(M). $$
By Definition  6.2  of $M\approx N$ and Theorem \ref{Thm6.6},
we have $\nu_{\om}(N)= \nu_{\om}(M)$, $S_N^-(\om)=S_M^-(\om)$, and
$S_N^+(\om)=S_M^+(\om)$. Thus in order to prove the right
inequalities of (\ref{eq2.81}) and (\ref{eq2.82}), it suffices to
prove them when $M$ is replaced by $N$. Then by Theorems 6.5 and 6.6
below, it suffices to prove them for each basic normal form given at
the beginning of Appendix 1 and listed in (\ref{eq6.10}) below.
This is a direct verification via Theorem 6.5 below and thus is left
to the readers. \enddemo

\specialnumber{2.2} \proclaim{Lemma}\lb{Lem2.2} Let $A\subset \U$ and $A'=\U\bs A${\rm .} Suppose $A$ is
symmetric with respect to $\R${\rm ;} i.e.{\rm ,} $\om\in A$ implies $\ol{\om}\in A$.
Then for any $M\in\Sp(2n)${\rm ,}  
\ee  \sum_{\om\in A}\left(\nu_{\om}(M) - S_M^-(\om)\right)
   + \sum_{\om\in A'}S_M^-(\om) \le \frac{e(M)}{2}.  \lb{eq2.83}\eee
\el

\demo{Proof} Fix $M\in\Sp(2n)$. Let $e(M,A)$ denote the total algebraic multiplicity
of all eigenvalues in $A\cap\sg(M)$. For any $\om\in\sg(M)\cap \U$, denote
its Krein type by $(p_{\om}(M), q_{\om}(M))$. Since $A$ is symmetric with
respect to $\R$, so is $A'$. Note that by the definition of the Krein type,
$$ p_{\om}(M) = q_{\ol{\om}}(M), \quad \hbox{for all }\om\in\U. $$
Therefore we have
\ee \sum_{\om\in A}p_{\om}(M) = \frac{e(M,A)}{2}, \quad
     \sum_{\om\in A'}q_{\om}(M) = \frac{e(M,A')}{2}. \lb{eq2.84}\eee
Together with Lemma \ref{Lem2.1}, this yields
$$ \sum_{\om\in A}\left(\nu_{\om}(M) - S_M^-(\om)\right)
    \le \frac{e(M,A)}{2}. $$
Using Lemma 4.10 of \cite{Lo7} and Theorem \ref{Thm6.6} below,
$$ \sum_{\om\in A'}S_M^-(\om) \le \frac{e(M,A')}{2}. $$
Thus (\ref{eq2.83}) holds. \enddemo

  Next we can give the following generalization of Theorem \ref{Thm2.2}.

\specialnumber{2.4}\proclaim{Theorem}\lb{Thm2.4} With $n\in\N${\rm ,} $\tau>0${\rm ,} and
$\ga\in\P_{\tau}(2n)${\rm ,} set  $M=\ga(\tau)${\rm .}  Extend $\ga$ to $[0,+\infty)$ by {\rm (\ref{eq1.18}).}
Then for any $m_1$ and $m_2\in\N${\rm ,} 

\aa
\noalign{\vskip-16pt}
&&\hskip-.25in \nu(\ga,m_1)+ \nu(\ga,m_2)- \nu(\ga, (m_1,m_2)) - \frac{e(M)}{2} \lb{eq2.90} \\
     &&\hskip.75in\le i(\ga,m_1+m_2) - i(\ga,m_1) - i(\ga,m_2)\nn \\
&&\hskip.75in\le \frac{e(M)}{2} + \nu(\ga, (m_1,m_2)) - \nu(\ga, m_1+m_2),
    \lb{eq2.91}\eaa
where $(m_1,m_2)$ is the greatest common divisor of $m_1$ and $m_2${\rm .}
\et

\demo{Proof} Fixing $m_1$ and $m_2\in\N$, define
\y
  \Psi_{m_1,m_2}(\th) &\hskip-5pt=\hskip-5pt& E\lrp{\frac{(m_1+m_2)\th}{2\pi}}
      - E\lrp{\frac{m_1\th}{2\pi}} - E\lrp{\frac{m_2\th}{2\pi}},
                \quad \hbox{for all }\th\in [0,2\pi), \\
 A &\hskip-5pt=\hskip-5pt& \lrc{\th\in [0,2\pi)\,|\,\frac{m_1\th}{2\pi}\in\Z
            \;{\rm or}\;\frac{m_2\th}{2\pi}\in\Z},     \\
 B &\hskip-5pt=\hskip-5pt& \lrc{\th\in [0,2\pi)\,|\, \frac{(m_1+m_2)\th}{2\pi}\in\Z
            \;{\rm and}\;\frac{m_1\th}{2\pi}\not\in\Z}, \\
 D &\hskip-5pt=\hskip-5pt& \lrc{\th\in [0,2\pi)\,|\,\Psi_{m_1,m_2}(\th)=0 }, \\
 D' &\hskip-5pt=\hskip-5pt& \lrc{\th\in (0,2\pi)\,|\,\Psi_{m_1,m_2}(\th)\not=0 }.
\ey
Note that we always have $\Psi_m(\th)=0$ or $-1$, $0\in A\subset D$
and $B\subset D'$.

  By Theorem \ref{Thm2.1}, we have
\aa
&& i(\ga,m_1+m_2) - i(\ga,m_1) - i(\ga,m_2) \lb{eq2.92}\\
&&\qquad\qquad = \lrp{S^+_M(1)+C(M)}
  + 2\sum_{\th\in(0,2\pi)}\Psi_{m_1,m_2}(\th)S^-_M\lrp{e^{\sqrt{-1}\th}} \nn\\
&&\qquad\qquad = S^-_M(1)+
 \left(\sum_{\th\in D\bs\{0\}}S^-_M\lrp{e^{\sqrt{-1}\th}}
 + \sum_{\th\in D'}S^-_M\lrp{e^{\sqrt{-1}\th}}\right) \nn\\
&&\qquad\qquad\quad  - 2\sum_{\th\in D'}S^-_M\lrp{e^{\sqrt{-1}\th}} \nn\\
&&\qquad\qquad = \sum_{\th\in D}S^-_M\lrp{e^{\sqrt{-1}\th}}
    - \sum_{\th\in D'}S^-_M\lrp{e^{\sqrt{-1}\th}}. \nn
\eaa

By Lemma \ref{Lem2.2}, noting that both $A$ and $B$ are symmetric
with respect to~$\R$, we obtain
\aa
\quad && \sum_{\th\in A}\nu_{e^{\sqrt{-1}\th}}(M)
  - \sum_{\th\in D}S_M^-\lrp{e^{\sqrt{-1}\th}}
  + \sum_{\th\in D'}S^-_M\lrp{e^{\sqrt{-1}\th}}  \lb{eq2.93}\\
&&\quad\qquad \le \sum_{\th\in A}\left(\nu_{e^{\sqrt{-1}\th}}(M)
       - S_M^-\lrp{e^{\sqrt{-1}\th}}\right)
 + \sum_{\th\in A'}S^-_M\lrp{e^{\sqrt{-1}\th}}\le \frac{e(M)}{2},
        \nn \\
&& \sum_{\th\in D}S_M^-\lrp{e^{\sqrt{-1}\th}}
 + \sum_{\th\in B}\nu_{e^{\sqrt{-1}\th}}(M)
 - \sum_{\th\in D'}S^-_M\lrp{e^{\sqrt{-1}\th}}  \lb{eq2.94}\\
&&\qquad\quad \le \sum_{\th\in B'}S_M^-\lrp{e^{\sqrt{-1}\th}}
 + \sum_{\th\in B}\left(\nu_{e^{\sqrt{-1}\th}}(M)
 - S^-_M\lrp{e^{\sqrt{-1}\th}}\right) \le \frac{e(M)}{2}.\nn
\eaa
By the definitions of $A$ and $B$, we have
\aa
\nu(\ga,m_1) + \nu(\ga,m_2) - \nu(\ga, (m_1,m_2))
  &=& \sum_{\th\in A}\nu_{e^{\sqrt{-1}\th}}(\ga),  \lb{eq2.95}\\
\nu(\ga,(m_1,m_2)) - \nu(\ga,m_1+m_2)
  &=& \sum_{\th\in B}\nu_{e^{\sqrt{-1}\th}}(\ga).  \lb{eq2.96}
\eaa
By (\ref{eq2.92}) we have
\aa
&&\quad i(\ga,m_1+m_2) - i(\ga,m_1) - i(\ga,m_2)   \lb{eq2.97}\\[4pt]
&&\hskip.75in = \sum_{\th\in D}S^-_M\lrp{e^{\sqrt{-1}\th}}
  - \sum_{\th\in D'}S^-_M\lrp{e^{\sqrt{-1}\th}}\nn \\
&&\hskip.75in = \sum_{\th\in A}\nu_{e^{\sqrt{-1}\th}}(M)
  \nn\\
&&\hskip.75in\quad - \left(\sum_{\th\in A}\nu_{e^{\sqrt{-1}\th}}(M)
        - \sum_{\th\in D}S^-_M\lrp{e^{\sqrt{-1}\th}}
      + \sum_{\th\in D'}S^-_M\lrp{e^{\sqrt{-1}\th}}\right). \nn
\eaa
Thus by (\ref{eq2.93}) and (\ref{eq2.96}), this implies (\ref{eq2.90}).

  Similarly, we have
\aa
&& i(\ga,m_1+m_2) - i(\ga,m_1) - i(\ga,m_2) \lb{eq2.98} \\
&&\qquad = \sum_{\th\in D}S^-_M\lrp{e^{\sqrt{-1}\th}}
  - \sum_{\th\in D'}S^-_M\lrp{e^{\sqrt{-1}\th}} 
  = \sum_{\th\in B}\nu_{e^{\sqrt{-1}\th}}(M)
 \nn\\
&&\qquad\quad + \left(\sum_{\th\in D}S_M^-\lrp{e^{\sqrt{-1}\th}}
 + \sum_{\th\in B}\nu_{e^{\sqrt{-1}\th}}(M)
 - \sum_{\th\in D'}S^-_M\lrp{e^{\sqrt{-1}\th}}\right). \nn
\eaa
Thus by (\ref{eq2.94}) and (\ref{eq2.96}), we obtain (\ref{eq2.91}).

  The proof is complete. \enddemo

\section{The injection map $p$}\lb{s3}

  In this section, we make precise the correspondence between effective
integers and iterations of  closed characteristics. This correspondence
is given by the map $p$ defined by the following lemma.

\specialnumber{3.1} \proclaim{Lemma}\lb{Lem3.1} Suppose $^\#\tilde\J(\Sigma,\alpha)< +\infty${\rm .} Then there
exist  an integer $K\ge 0$ and an injection map
$p: \N+K \to \V_{\infty}(\Sg,\alpha)\times\N$ such that
\bd
\item[{\rm (i)}] For any $k\in\N+K${\rm ,} $(\tau,x)\in\J(\Sigma,\alpha)$ and $m\in\N$
satisfying $p(k)=([(\tau,x)],m)${\rm , (\ref{eq1.29})} and {\rm (\ref{eq1.30})} hold{\rm ,} and
\item[{\rm (ii)}] For any $k_j\in\N+K${\rm ,} $k_1<k_2${\rm ,} $(\tau_j,x_j)\in\J(\Sg,\alpha)$ satisfying
$p(k_j)=([(\tau_j,x_j)],m_j)$ with $j=1,2${\rm ,}
 \vglue4pt
\hglue-28pt {\rm (3.1)} \hfil  ${\displaystyle \hat{i}(x_1,m_1) < \hat{i}(x_2,m_2). }$\hfill
\ebd
\el

\advance\eqcount by 1
\demo{Proof}  (i)  By Theorem V.3.4 of  \cite{Ek3}, for each $k\in\N$, there is a
$(\tau,x)\in\J(\Sg,\alpha)$ such that $(\tau,x)$ is $(m,k)$-variationally
visible for some $m\in\N$, i.e. $[(\tau,x)]\in\V(\Sigma,\alpha)$. We define
a map $p_1: \N\to\V(\Sigma,\alpha)\times\N$  by  $p_1(k)=([(\tau,x)],m)$.
Since $\td{\J}(\Sigma)< +\infty$, (\ref{eq1.38}) holds by  Theorem V.3.4
of  \cite{Ek3}, i.e. $c_j<c_k<0$ whenever $j<k$. Thus if
$p_1(j)=p_1(k)=([\tau,x)],m)$ for some $j<k$, by  (\ref{eq1.29}), we obtain
$$  c_j = f(u_m^x) = c_k.  $$
This contradiction proves that $p_1$ is injective.

  Since $^{\#}\td{\J}(\Sg)< +\infty$, there exists an integer  $K\ge 0$ such that
all critical values $c_{k+K}$ with $k\in\N$ come from iterations of  elements
in  $\V_{\infty}(\Sigma,\alpha)$. Thus $p_1(k+K)\in\V_{\infty}(\Sg,\alpha)\times\N$
for any $k\in\N$. We define
\ee  p(k)=p_1(k),  \qquad \hbox{for all }k\in\N+K.  \lb{eq3.1a}\eee
Then $p$ is injective, and (\ref{eq1.29}) and (\ref{eq1.30}) hold.

\vglue4pt
 (ii) By Theorem V.3.11 of \cite{Ek3} we have
$\ga\equiv\ga^-(\Sigma)=\ga^+(\Sigma)>0$, where $\ga^{\pm}(\Sigma)$ are defined
by (V.3.62)--(V.3.65) of \cite{Ek3} as follows:
\aa
\ga^+(\Sigma) &=& C_{\alpha}^{-1}\limsup_{k\to\infty}((-c_k)^{(2-\alpha)/\alpha}k)^{-1},  \nn\\
 \ga^-(\Sigma) &=& C_{\alpha}^{-1}\liminf_{k\to\infty}((-c_k)^{(2-\alpha)/\alpha}k)^{-1},  \nn
\eaa
with $C_{\alpha}=\frac{4}{\alpha}(1-\frac{\alpha}{2})^{\frac{\alpha-2}{\alpha}}$. Thus we have
\ee  \lim_{k\to +\infty} k |c_k|^{\frac{2-\alpha}{\alpha}}C_{\alpha}
     = \frac{1}{\gamma}.  \lb{eq3.1b}\eee

  For $k\in\N+K$, let $p(k)=([(\tau,x)],m)$ for some
$[(x,\tau)]\in\V_{\infty}(\Sigma,\alpha)$ and $m\in\N$. Then $(x,\tau)$ is
$(m,k)$-variationally visible. So by (V.3.45) of \cite{Ek3},
\aa
 c_k &=& f(u_m^x),  \nn\\
   |c_k|^{\frac{2-\alpha}{\alpha}} &=& 2(C_{\alpha}mA(\tau,x))^{-1}
         =2(C_{\alpha}A(m\tau,x^m))^{-1},  \nn\eaa
where $A(\tau,x)=\frac{1}{2}\int_0^{\tau}(\dot{x},Jx)dt$. Note that by
Lemma V.3.12 of  \cite{Ek3} and (\ref{eq1.19a}), there holds
$$ \frac{\hat{i}(x,m)}{A(m\tau,x^m)} = \frac{\hat{i}(x,1)}{A(\tau,x)}
=\frac{1}{\ga}.  $$
This implies
$$ |c_k|^{\frac{2-\alpha}{\alpha}} \hat{i}(x,m) = \frac{2}{\ga C_{\alpha}}. $$
Since $c_k<c_{k+1}<0$ for $k\in\N$, we get our results.
\enddemo
\pagebreak

  By  Theorem \ref{Thm2.3} and Corollary \ref{Cor1.2}  we have

\specialnumber{3.1}\proclaim{{C}orollary}\lb{Cor3.1} Fix $\Sigma\in\H(2n)$ and $\alpha\in (1,2)${\rm .}
For any $(\tau,x)\in\J(\Sigma,\alpha)$ and $m\in\N${\rm ,} 
\aa
i(x,m+1) - i(x,m) &\ge& 2,  \lb{eq3.2}\\
\quad \qquad i(x,m+1) + \nu(x,m+1) - 1 &\ge&  i(x,m+1) > i(x,m) + \nu(x,m) - 1.
    \lb{eq3.3}\\
\hat{i}(x,1) &\ge& 2. \lb{eq3.4}\eaa
\ec

\demo{Proof} It suffices to prove (\ref{eq3.4}). By (\ref{eq3.2}) and an induction argument,
we obtain
$$ \frac{i(x,m)}{m}\ge \frac{i(x,1)+2m-2}{m}, \qquad \hbox{for all }m\in\N.  $$
This implies (\ref{eq3.4}). \enddemo

\demo{{R}emark {\rm 3.1}}  The sharpest estimate $\hat{i}(x,1) > 2$ for every
$(\tau,x)\in\J(\Sg,\alpha)$ was first proved in Theorem 2 of \cite{EH} (cf.\ also
Theorem I.7.7 of \cite{Ek3}). A different proof was given by Lemma 6.7 of
\cite{Lo9}. A third proof can be given by use of     Corollary \ref{Cor2.1} above
and   the study of splitting numbers in \cite{Lo7}.  \er

\vglue-6pt

\section{The common index jump of closed characteristics}\lb{s4}
\vglue-6pt

  The goal of this section is to prove the common index jump claim (\ref{eq1.44})
as well as other related results.

\demo{{\rm 4.1.} A common selection theorem} \enddemo

\specialnumber{4.1}\proclaim{Theorem}\lb{Thm4.1} Fix an integer $q>0$. Let $\mu_i\ge 0$ and $\beta_i$ be integers
for all $i=1,\ldots,q$. Let $\alpha_{i,j}$ be positive numbers for $j=1,\ldots,\mu_i$
and $i=1,\ldots,q$. Let $\delta\in(0,\frac{1}{2})$ satisfying
\ee  \dl \max_{1\le i\le q}\mu_i  < \frac{1}{2}. \lb{eq4.1a}\eee
Set
\ee D_i=\beta_i+\sum_{j=1}^{\mu_i}\alpha_{i,j}, \qquad
\ox{{\it for}}\;i=1,\ldots,q.  \lb{eq4.1}\eee
Suppose
\ee  D_i>0,\qquad \hbox{for all }i=1,\ldots,q. \lb{eq4.2}\eee
Then there exist infinitely many $(N, m_1,\ldots,m_q)\in\N^{q+1}$ such that
\aa
m_i\beta_i+\sum_{j=1}^{\mu_i}E(m_i\alpha_{i,j})  = N + \Delta_i, & &
         \hbox{for all }i=1,\ldots, q,    \lb{eq4.4}\\
\qquad\quad\min \{\{m_i\alpha_{i,j}\}, \;1-\{m_i\alpha_{i,j}\}\} < \delta, & &
         \hbox{for all }j=1,\ldots,\mu_i, i=1,\ldots,q,   \lb{eq4.5}\\
m_i\alpha_{i,j}\in\N & &\ox{if}\;\alpha_{i,j}\in\Q,   \lb{eq4.6}\eaa
where
\ee \Delta_i = \sum_{0<\{m_i\alpha_{i,j}\}<\delta}1 \lb{eq4.7}\eee
for all $i=1,\ldots, q$.  \et

\demo{{P}roof} Firstly, we reduce the claims (\ref{eq4.4})--(\ref{eq4.6}) to a
dynamical problem on a torus. Then we prove the existence of numbers
$(N, m_1,\ldots,m_q)$ by using properties of closed additive subgroups
of tori. The proof is carried out in two steps.
\enddemo

\demo{Step 1} {\it Reduction to a problem on the  torus}.
  We consider the left-hand side of (\ref{eq4.4}) first. Using the function $\phi(\cdot)$ defined
at the beginning of Section~\ref{s2}, we obtain
\begin{equation}
m_i\beta_i + \sum_{j=1}^{\mu_i}E(m_i\alpha_{i,j}) 
= m_iD_i + \sum_{j=1}^{\mu_i}\left(\phi(m_i\alpha_{i,j}) - \{m_i\alpha_{i,j}\}\right).
 \lb{eq4.11} 
\end{equation}

  To handle all the rational rotators, for $1\le i\le q$ we require each $m_i$ having a factor
$M\in\N$ such that $M\alpha_{i,j}\in\N$ whenever  $\alpha_{i,j}\in\Q$ for $j=1,\ldots,\mu_i$
and $i=1,\ldots,q$. Let $M=1$ if  no such $\alpha_{i,j}$ exists.

  To get the common integer $N\in\N$, we replace $m_iD_i$ in (\ref{eq4.11})
by $\frac{N}{MD_i}MD_i$. To make $m_i$ an integer, we subtract the decimal
part of $\frac{N}{MD_i}$ from itself. Because of the requirement of our
later torus problem and to keep more flexibility in our choice of $N$, we add
a term $\chi_i\in \{0,1\}$ to it. This yields
\aa
  m_iD_i
&=& \frac{N}{MD_i}MD_i - \lrc{\frac{N}{MD_i}}MD_i + \chi_iMD_i\lb{eq4.12}\\
&=& N + \lrp{\chi_i - \lrc{\frac{N}{MD_i}MD_i}}MD_i.  \nn\eaa
Here we define
\ee  m_i = \left(\lrb{\frac{N}{MD_i}}
+\chi_i\right)M, \qquad \hbox{for all }i=1, \ldots, q, \lb{eq4.14}\eee
where $\chi_i=0$ or $1$ for $1\le i\le q$ will be determined later. By this choice, (\ref{eq4.11})
becomes
\aa 
m_i\beta_i+ \sum_{j=1}^{\mu_i}E(m_i\alpha_{i,j})  
&=& N + \lrp{\chi_i - \lrc{\frac{N}{MD_i}MD_i}}MD_i\lb{eq4.16} \\
 &&      + \sum_{j=1}^{\mu_i}\left(\phi(m_i\alpha_{i,j}) - \{m_i\alpha_{i,j}\}\right).
   \nn\eaa

  Now the claims (\ref{eq4.5}) and (\ref{eq4.6}) require that only the following
three possibilities for each term $\{m_i\alpha_{i,j}\}$ can happen:\pagebreak

  (A) $\;\{m_i\alpha_{i,j}\}=0$ if $\alpha_{i,j} \in\Q$,

\vglue6pt

  (B) $\;0 < \{m_i\alpha_{i,j}\}<\dl$ if $\alpha_{i,j}\in \R\bs\Q$,  or
\vglue6pt

  (C) $\;1-\dl < \{m_i\alpha_{i,j}\} < 1$  if $\alpha_{i,j}\in \R\bs\Q$.
\vglue6pt

  Note that (A) already holds by our choice of $m_i$'s of (\ref{eq4.14}).
Suppose now that these three requirements are fulfilled by our choice of $m_i$'s.
Then we have
\aa &&\lb{eq4.18}\\
&& \hskip-24pt \sum_{j=1}^{\mu_i}\left(\phi(m_i\alpha_{i,j}) - \{m_i\alpha_{i,j}\}\right)\nn\\
&& \quad = \sum_{A} + \sum_{B}
       + \sum_{C} \left(\phi(m_i\alpha_{i,j}) - \{m_i\alpha_{i,j}\}\right) \nn\\
&& \quad = \sum_{0<\{m_i\alpha_{i,j}\}<\dl}1 -  \sum_{0<\{m_i\alpha_{i,j}\}<\dl}\{m_i\alpha_{i,j}\}
                   + \sum_{0< 1-\{m_i\alpha_{i,j}\}<\dl}(1- \{m_i\alpha_{i,j}\})  \nn\\
&& \quad = \Delta_i - \sum_{0<\{m_i\alpha_{i,j}\}<\dl}\{m_i\alpha_{i,j}\}
                   + \sum_{0< 1-\{m_i\alpha_{i,j}\}<\dl}(1- \{m_i\alpha_{i,j}\}).  \nn\eaa
Thus (\ref{eq4.16}) becomes
\aa\qquad\quad
m_i\beta_i &+& \sum_{j=1}^{\mu_i}E(m_i\alpha_{i,j})  \lb{eq4.20}\\
&=& N + \lrp{\chi_i - \lrc{\frac{N}{MD_i}MD_i}}MD_i   \nn\\
& &\;   +  \Delta_i - \sum_{0<\{m_i\alpha_{i,j}\}<\dl}\{m_i\alpha_{i,j}\}
                   + \sum_{0< 1-\{m_i\alpha_{i,j}\}<\dl}(1- \{m_i\alpha_{i,j}\}).  \nn\eaa
Therefore
\ee
  \left| m_i\beta_i + \sum_{j=1}^{\mu_i}E(m_i\alpha_{i,j}) - N - \Delta_i\right|
    \le \left| \chi_i - \lrc{\frac{N}{MD_i}}\right|MD_i + \mu_i\dl.  \lb{eq4.22}\eee
Here, to get (\ref{eq4.4})  by (\ref{eq4.1a})  we need
\ee  \left| \chi_i - \lrc{\frac{N}{MD_i}}\right|MD_i   < \frac{1}{2}.  \lb{eq4.24}\eee

  Next we estimate $\{m_i\alpha_{i,j}\}$. By our choice of $m_i$ in (\ref{eq4.14}),  
\aa
\{m_i\alpha_{i,j}\}
&=& \lrc{M\lrp{\lrb{\frac{N}{MD_i}}+\chi_i}\alpha_{i,j}}   \lb{eq4.26}\\[4pt]
&=& \lrc{\frac{N\alpha_{i,j}}{D_i} + \lrp{\chi_i - \lrc{\frac{N}{MD_i}}}M\alpha_{i,j}}  \nn\\[4pt]
&=&\{A_{i,j}(N)+B_{i,j}(N)\},   \nn\eaa
\pagebreak

\noindent
where in (\ref{eq4.26}), we set
\ee    A_{i,j}(N) = \lrc{\frac{N\alpha_{i,j}}{D_i}} - \chi_{i,j},\qquad
        B_{i,j}(N) = \lrp{\chi_i - \lrc{\frac{N}{MD_i}}}M\alpha_{i,j}, \hskip.25in \lb{eq4.28}\eee
and $\chi_{i,j}=0$ or $1$ will be determined later.

  Assume now that we can choose $N\in\N$ such that
\aa
\left| \lrc{\frac{N\alpha_{i,j}}{D_i}} - \chi_{i,j}\right|
        &=& |A_{i,j}(N)| < \frac{\dl_1}{3},  \lb{eq4.30}\\
\left| \lrp{\chi_i - \lrc{\frac{N}{MD_i}}}M\alpha_{i,j}\right|
         &=& |B_{i,j}(N)| < \frac{\dl_1}{3},  \lb{eq4.32}\eaa
for a given $\dl_1$ satisfying $0<\dl_1<\dl<1/2$.  Then if $A_{i,j}(N)+B_{i,j}(N) \ge 0$,
we obtain
$$ \{m_i\alpha_{i,j}\}=A_{i,j}(N)+B_{i,j}(N) <\dl_1 <\dl.  $$
If  $A_{i,j}(N)+B_{i,j}(N) < 0$, we obtain
$$ \{m_i\alpha_{i,j}\}=A_{i,j}(N)+B_{i,j}(N) - (-1). $$
That is,
$$ 0 < 1 - \{m_i\alpha_{i,j}\}= - (A_{i,j}(N)+B_{i,j}(N)) <\dl_1 <\dl.  $$

  Therefore to prove (\ref{eq4.4})--(\ref{eq4.6}), by (\ref{eq4.24}), (\ref{eq4.30}), and (\ref{eq4.32}),
it suffices to prove that we can choose integers $\chi_i$ and $\chi_{i,j}$ to be $0$ or $1$ and
choose infinitely many integers $N\in\N$ such that all the quantities
\ee  \left| \lrc{\frac{N\alpha_{i,j}}{D_i}} - \chi_{i,j}\right|  \quad {\rm and}\quad
        \left| \lrc{\frac{N}{MD_i}} - \chi_i \right|   \lb{eq4.34}\eee
can be made simultaneously as small as we want.

  Let $n=q+\sum_{i=1}^q\mu_i$, and
\ee  v = (\frac{1}{MD_1}, \ldots, \frac{1}{MD_q}, \frac{\alpha_{1,1}}{D_1},
          \frac{\alpha_{1,2}}{D_1}, \ldots, \frac{\alpha_{1,\mu_1}}{D_1},
          \frac{\alpha_{2,1}}{D_2}, \ldots, \frac{\alpha_{q,\mu_q}}{D_q}) \in \R^n. \lb{eq4.36}\eee
Then the problem becomes, for any given small $\ep \in (0,\min\{\dl,1/3\})$, to find a vertex $\chi$
of the cube $[0,1]^n$ and infinitely many integers $N\in \N$ such that 
\ee   | \{Nv\}-\chi | < \ep.  \lb{eq4.38}\eee
This is a problem of the dynamics on the standard torus $\T^n=\R^n/\Z^n$.
\enddemo

\demo{Step 2} {\it Dynamics on a torus}.
  To solve the problem (\ref{eq4.38}), it suffices to note that the closure
$\ol{G}$ of the set $G=\{\{mv\}\,|\,m\in\Z\}$ in $\T^n$ forms a closed additive
subgroup of the torus $\T^n$. Thus $\ol{G}$ is a product of a possibly lower
dimensional torus and a cyclic group, and $G$ contains the identity element of~$\T^n$.
\pagebreak

  Define a dynamical system $f:\T^n\to\T^n$ by $f(x)=x+\pi(v)$ for all
$x\in\T^n$. We consider the $\alpha$ and $\omega$-limit sets of $f$ which
are defined as usual by
$$ \alpha(x)=\bigcap_{N\in\N}\ol{\{ f^k(x)\,|\,k\le -N\}}\quad{\rm and}
    \quad\omega(x)=\bigcap_{N\in\N}\ol{\{ f^k(x)\,|\, k\ge N\}},
 $$
{for all } $x\in\T^n. $

  If $v\in\Q^n$, $G$ must be the finite cyclic group generated by $\pi(v)$. So
we have $\alpha(0)=\omega(0)=G=\ol{G}$.

  If $v\in \R^n\bs\Q^n$, then we have $^{\#}G= +\infty$ and hence the limit
sets $\alpha(0)$ and $\omega(0)$ are nonempty. By the definition of the limit
sets and the fact $G$ is an additive subgroup of $\T^n$, we have
$$ \alpha(0)+G\subset\alpha(0),\qquad \omega(0)+G\subset\omega(0).$$
Since the limit sets are closed,  
$$ \alpha(0)+\ol{G}\subset\alpha(0),\qquad \omega(0)+\ol{G}\subset\omega(0).$$
Thus by the facts that $\alpha(0)\cup\omega(0)\subset \ol{G}$ and that $\ol{G}$
is a closed additive subgroup of $\T^n$,  
\ee \alpha(0)=\omega(0)=\ol{G}. \lb{eq4.39}\eee

  Therefore there always exist infinitely many $N\in\N$ such that the point
$\{Nv\}$ is located in the open ball in $\T^n$ centered at its identity element
with radius $\ep<1/3$. More precisely, we have proved that in the product
space $X=[0,1]^n$, there exist one vertex $x$ of $X$ and infinitely many
$N\in \N$ such that
\ee  \{Nv\}\; \in \; B_{\ep}(x)\cap X,  \lb{eq4.40}\eee
where $B_{\ep}(x)=\{y\in\R^n\,|\, |y-x| <\ep\}$. Now we define $\chi=x$, and then
(\ref{eq4.38}) holds.

  The proof is complete. \hfill\qed\enddemo

  In order to prove Theorem \ref{Thm1.4}, we need to know more about the possible
choices of the vector $\chi$ in Theorem \ref{Thm4.1}. The following result gives this
information.

\specialnumber{4.2}\proclaim{Theorem}\lb{Thm4.2} Fix $v=(v_1,\ldots,v_n)\in \R^n${\rm .} Let $H$ be the
closure of
$\{ \{mv\} \,|\, m\in\N\}$ in $\T^n$ and $V=T_0\pi^{-1}H$ be the tangent space of
$\pi^{-1}H$ at the origin in $\R^n${\rm ,} where $\pi:\R^n\to \T^n$ is the projection map{\rm .}
Define
\ee  A(v) =V\bs\mathbold{\cup}_{v_k\in\R\bs\Q}\{x=(x_1,\ldots,x_n)\in V\,|\, x_k=0\}.
\lb{eq4.42}\eee
Define $\psi(x)=0$ when $x\ge 0$ and $\psi(x)=1$ when $x<0${\rm .} Then for \pagebreak
any $a=(a_1,\ldots,a_n) \in A(v)${\rm ,} the vector
\ee  \chi = (\psi(a_1), \ldots, \psi(a_n)) \lb{eq4.44}\eee
makes {\rm (\ref{eq4.38})} holds for infinitely many $N\in\N${\rm .}

\def\ritem#1{\item[{\rm #1}]}

  Moreover{\rm ,} this set $A(v)$ possesses the following properties\/{\rm :}
\begin{itemize}
\ritem{(a)} $A(v)\ne\emptyset${\rm .}

\ritem{(b)} When $v\in \Q^n${\rm ,}   then $V=A(v)=\{0\}${\rm .}

\ritem{(c)} When $v\in\R^n\bs\Q^n${\rm ,}  then $\dim V\geq 1$, $0\not\in A(v)\subset V${\rm ,}
$A(v)=-A(v)${\rm ,} and   $A(v)$ is open in $V${\rm .}

\ritem{(d)} When $\dim V=1${\rm ,} then  $A(v)=V\bs\{0\}${\rm .}

\ritem{(e)} When $\dim V\ge 2${\rm ,} $A(v)$ is obtained from $V$ by deleting all the
coordinate hyperplanes with dimension strictly smaller than $\dim V$ from
$V${\rm ,} especially $\dim A(v)=\dim V${\rm .}
\end{itemize}

\et

\demo{Proof} If $v\in\Q^n$, $H$ is a cyclic subgroup of $\T^n$. We then have
$V=\{0\}$. Thus $A(v)=V=\{0\}$.

  Suppose $v=(v_1,\ldots,v_n)\in\R^n\bs\Q^n$. Then $H$ is a closed additive
subgroup of $\T^n$ which is a product of a torus with a cyclic subgroup of
$\T^n$. By the definition of $V$, we have $\dim V\ge 1$. For any
$a=(a_1,\ldots,a_n)\in V$,
$$ a_i \ne  0 \quad {\rm only\;if\;} v_i \in \R\bs \Q. $$
When $v_i \in \R\bs \Q$,  there exists $a=(a_1,\ldots,a_n)\in V$ with
$a_i\not= 0$. Thus for any $a\in A(v)$, the point $\chi$ defined by
(\ref{eq4.44}) gives a vertex $\chi$ of $X=[0,1]^n$ which is in the closure of
$\{\{mv\}\,|\,m\in\N\}$ in $X$. Therefore (\ref{eq4.38}) holds for infinitely
many $N\in\N$ and for this $\chi$.

  The claims (a) to (e) follow from this argument, and their proofs are therefore
omitted.  \enddemo

\demo{{\rm 4.2.} The common index jump of symplectic paths} 
  The following is the main result of this subsection.
\enddemo

\specialnumber{4.3}\proclaim{Theorem}\lb{Thm4.3} Let $\ga_k\in\P_{\tau_k}(2n)$ for $k=1,\ldots,q$ be
a finite collection of symplectic paths{\rm .} Let $M_k=\ga(\tau_k)${\rm .} Extend $\ga_k$ to
$[0,+\infty)$ by iteration via {\rm (\ref{eq1.18})} for $k=1,\ldots,q${\rm .} Suppose
\ee  \hat {i}(\ga_k,1)>0, \qquad \hbox{for all }k=1,\ldots,q.  \lb{eq4.50}\eee
Then there exist infinitely many $(N, m_1,\ldots,m_q)\in\N^{q+1}$ such that
\aa
   \nu(\ga_k,2m_k-1)  &=& \nu(\ga_k,1),    \lb{eq4.51}\\
   \nu(\ga_k,2m_k+1) &=& \nu(\ga_k,1),    \lb{eq4.52}\\
&&\hskip-1.75in i(\ga_k,2m_k-1) + \nu(\ga_k,2m_k-1)\lb{eq4.53}\\
 &=& 2N
        - \left(i(\ga_k,1) + 2S^+_{M_k}(1) - \nu(\ga_k,1)\right),        \nn\\
   i(\ga_k,2m_k+1) &=& 2N + i(\ga_k,1),   \lb{eq4.54}\\
   i(\ga_k,2m_k) &\ge& 2N-\frac{e(M_k)}{2}\ge 2N-n,   \lb{eq4.55}\\
 \qquad  i(\ga_k,2m_k) + \nu(\ga_k,2m_k) &\le& 2N+\frac{e(M_k)}{2}\le 2N+n,   \lb{eq4.56}\eaa
for every $k=1,\ldots,q$.
\et

{\it Proof}. We complete the proof in four steps.
\demo{Step 1}  {\it Application of Theorem} \ref{Thm4.1}.  
  Set
\aa
\delta_0 &=& \min_{1\le k\le q}\lrc{\frac{1}{2}, \frac{\theta}{2\pi},1-\frac{\theta}{2\pi}
    \mid\theta\in(0,2\pi)\;\ox{{\it and}}\;e^{\sqrt{-1}\theta}\in\sigma(M_k)},  \lb{eq4.56a}\\
C_k &=& \sum_{\theta\in(0,2\pi)}S^-_{M_k}\lrp{e^{\sqrt{-1}\theta}},  \lb{eq4.56b}\\
\rho_k  &=& i(\ga_k,1) +S^+_{M_k}(1)-C_k,   \lb{eq4.56c}\\
\qquad \quad I(k,m) &=& m\rho_k +
\sum_{\theta\in(0,2\pi)}E\lrp{\frac{m\theta}{\pi}}S^-_{M_k}\lrp{e^{\sqrt{-1}\theta}}
    \lb{eq4.56d}\eaa
for $k=1,\ldots,q$ and $m\in\N$, where $\sigma(M_k)$ denotes the spectrum of $M_k$.

  By the definition of $I(k,m)$, we rewrite it as
\ee I(k,m) = m\rho_k+\sum_{\theta\in(0,2\pi)}
\sum_{j=1}^{S^-_{M_k}\lrp{e^{\sqrt{-1}\theta}}}E\lrp{\frac{m\theta}{\pi}}.
                \lb{eq4.57} \eee
Note that by (\ref{eq4.50}) and Corollary \ref{Cor2.1}, for every $k=1,\ldots,q$,
\ee  0<\hat{i}(\ga_k,1)=\rho_k
  +\sum_{\theta\in(0,2\pi)}\frac{\theta}{\pi}S^-_{M_k}\lrp{e^{\sqrt{-1}\theta}}.  \lb{eq4.58}\eee

  Now in the statement of  Theorem \ref{Thm4.1}, other than the integer $q$, we set
\aa
 \delta &\in& (0,\delta_0), \quad \beta_i = \rho_i,
      \lb{eq4.58a} \\
 \mu_i& =& \sum_{\theta\in (0,2\pi)}S^-_{M_i}\lrp{e^{\sqrt{-1}\theta}},
       \quad       D_i=\hat{i}(\ga_i,1),  \quad            \hbox{for all }1\le i\le q,   \nn\\
\qquad\quad  \alpha_{i,j} &=& \frac{\theta_j}{\pi},  \quad \mbox{where} \;
            e^{\sqrt{-1}\theta_j}\in\sigma(M_i),\quad  \hbox{for all }1\le j\le \mu_i,\, 1\le i\le q.
                    \lb{eq4.58b}\eaa
Note that $\mu_i$ is a nonnegative integer by Corollary 4.13 of \cite{Lo7}.\pagebreak

  By (\ref{eq4.58}), the condition (\ref{eq4.2}) of Theorem \ref{Thm4.1} holds.
Applying this theorem to $I(k,m)$'s for $k=1,\ldots,q$ and any
$\delta\in (0,\delta_0)$, there exist infinitely many $(N,m_1,\ldots,m_q)\in\N^{q+1}$
such that
\aa
&&  I(k,m_k) = N + \Delta_k,    \lb{eq4.59}\\
&& \min\lrc{\lrc{\frac{m_k\theta}{\pi}}, 1-\lrc{\frac{m_k\theta}{\pi}}}
       <  \delta, \qquad  \mbox{if}\;  e^{\sqrt{-1}\theta}\in\sigma(M), \lb{eq4.60}\\
&& \frac{m_k\theta}{\pi} \in \Z,  \qquad
            \mbox{if}\;\;\frac{\theta}{\pi}\in\Q\cap(0,2)\;\;\mbox{and}\;\;
                     e^{\sqrt{-1}\theta}\in\sigma(M_k),     \lb{eq4.61}\\
&&\Delta_k=\sum_{0<\{m_k\theta/\pi\}<\delta}S^-_{M_k}\lrp{e^{\sqrt{-1}\theta}},
      \lb{eq4.62}\eaa
for $k=1,\ldots, q$.
\enddemo

\demo{Step 2} {\it Verifications of}\/ (\ref{eq4.51}) {\it and} (\ref{eq4.52}).  
  Whenever $e^{\sqrt{-1}\theta}\in\sigma(M_k)$ and $\frac{\th}{\pi}\in \Q\cap (0,2)$,
by (\ref{eq4.61}), we always have $2m_k\theta \in 2\pi\Z$. Thus
for any such~$\theta$,
$$ 2m_k\theta \pm \theta\not\in 2\pi\Z. $$
Since the change of the nullity happens only
when iterations of some eigenvalues in $\U\bs\{1\}$ hit $1$, we
obtain
$$ \nu(\ga_k,2m_k-1) = \nu(\ga_k,2m_k+1) = \nu(\ga_k,1).   $$
Thus (\ref{eq4.51}) and (\ref{eq4.52}) hold.
\enddemo

\demo{Step 3} {\it Verifications of} (\ref{eq4.53}) {\it and} (\ref{eq4.54}).  
  By Theorem \ref{Thm2.1}, (\ref{eq4.59}), (\ref{eq4.51}), and the definition of $I(k,m)$,
\aa
&&  i(\ga_k,2m_k-1) + \nu(\ga_k,2m_k-1)  \lb{eq4.63} 
\\
&&\qquad\quad= 2I(k,m_k) - \lrp{i(\ga_k,1) + S^+_{M_k}(1) - C_k} - \lrp{S^+_{M_k}(1) + C_k}  \nn\\
& & \qquad\qquad  - 2\sum_{\theta\in(0,2\pi)}
              \xi_-(m_k,\theta)S^-_{M_k}\lrp{e^{\sqrt{-1}\theta}} +  \nu(\ga_k,1)   \nn\\
&&\qquad\quad =2(N + \Delta_k) -\lrp{i(\ga_k,1) + 2S^+_{M_k}(1)}  \nn\\
& & \qquad\qquad  -2\sum_{\theta\in(0,2\pi)}
           \xi_-(m_k,\theta)S^-_{M_k}\lrp{e^{\sqrt{-1}\theta}} + \nu(\ga_k,1), \nn \eaa
where we define
\aa
&&\lb{eq4.64} \\
\xi_-(m_k,\theta)
&=& E\lrp{\frac{2m_k\theta}{2\pi}} - E\lrp{\frac{(2m_k-1)\theta}{2\pi}}  \nn\\
&=& E\lrp{\lrb{\frac{m_k\theta}{\pi}} + \lrc{\frac{m_k\theta}{\pi}}} -
          E\lrp{\lrb{\frac{m_k\theta}{\pi}} + \lrc{\frac{m_k\theta}{\pi}} - \frac{\theta}{2\pi}}  \nn\\
&=& E\lrp{\lrc{\frac{m_k\theta}{\pi}}} - E\lrp{\lrc{\frac{m_k\theta}{\pi}} - \frac{\theta}{2\pi}}.
         \nn\eaa
Note that $\xi_-(m_k,\th)$ takes only the value $0$ or $1$. To evaluate $\xi_-(m_k,\theta)$
for any $\theta\in(0,2\pi)$ and $e^{\sqrt{-1}\theta}\in\sigma(M_k)$ with some
$k\in \{1,\ldots,q\}$,  we consider the following three cases:
\enddemo 

\demo{Case 1} $\{\frac{m_k\theta}{\pi}\}=0$.
  By (\ref{eq4.64}) and the definition of the function $E(\cdot)$, we have
$\xi_-(m_k,\theta)=0$.
\enddemo

\demo{Case 2} $\{\frac{m_k\theta}{\pi}\}\in (0,\delta)$.
  Since $\delta\in (0,\delta_0)$, by the definition (\ref{eq4.56a}) of $\delta_0$, we obtain
$$   \lrc{\frac{m_k\theta}{\pi}} - \frac{\theta}{2\pi} < \delta_0  - \frac{\theta}{2\pi} \le 0. $$
Thus, $\xi_-(m_k,\theta)=1$.
\enddemo

\demo{Case 3} {\it $\{\frac{m_k\theta}{\pi}\}\ge \delta$ and} $S^-_{M_k}\lrp{ e^{\sqrt{-1}\theta}}>0$.
  By  (\ref{eq4.60}),   
$$ 1-\delta < \lrc{\frac{m_k\theta}{\pi}} < 1  $$
must hold.
Thus by the fact $\delta\in (0,\delta_0)$ and the definition (\ref{eq4.56a}) of $\delta_0$,  
$$  \lrc{\frac{m_k\theta}{\pi}} - \frac{\theta}{2\pi} > 1 - \delta - \frac{\theta}{2\pi} > 0.  $$
So in this case we have $\xi_-(m_k,\theta)=0$.

  Therefore only in the above case 2, the term $\xi_-(m_k,\th)$ makes contribution. Together
with (\ref{eq4.63}) and (\ref{eq4.64}),
\aa
 i(\ga_k,2m_k-1)  +  \nu(\ga_k,2m_k-1)   
&=&  2(N + \Delta_k) - \lrp{i(\ga_k,1) + 2S^+_{M_k}(1)} \nn\\
  & & - \ 2\Delta_k + \nu(\ga_k,1)   \nn\\
&=&  2N - \lrp{i(\ga_k,1) + 2S^+_{M_k}(1) - \nu(\ga_k,1) }.  \nn\eaa
This proves (\ref{eq4.53}).
\enddemo

  Similarly,
\aa &&\lb{eq4.65}\\
i(\ga_k,2m_k+1)
&\hskip-3pt =\hskip-3pt&  2I(k,m_k) + \lrp{i(\ga_k,1) + S^+_{M_k}(1) - C_k}\nn\\
&& -\ \lrp{S^+_{M_k}(1) +
C_k}    
   + 
2\sum_{\th\in(0,2\pi)}\xi_+(m_k,\th)S^-_{M_k}\lrp{e^{\sqrt{-1}\theta}},
\nn\eaa
where $\xi_+(m_k,\th)= E(\frac{(2m_k+1)\theta}{2\pi}) - E(\frac{2m_k\theta}{2\pi})$.
Similar  to our  discussion above,  $\xi_+(m_k,\th)=1$ in    cases 1 and 3 above,
and $\xi_+(m_k,\th)=0$ in   case 2. Thus from (\ref{eq4.65}) we obtain
\aa
\quad i(\ga_k,2m_k+1)
&=&  2(N + \Delta_k) + (i(\ga_k,1) - 2C_k) + 2(C_k - \Delta_k)\lb{eq4.66}\\
&=&  2N + i(\ga_k,1).   \nn\eaa
This proves (\ref{eq4.54}).

\demo{Step 4} {\it Verifications of} (\ref{eq4.55}) {\it and} (\ref{eq4.56}). 
  By  the inequality (\ref{eq2.21}) in Theorem \ref{Thm2.2} and  (\ref{eq4.52})
as well as  (\ref{eq4.54}),   
\aa
i(\ga_k,2m_k)
&\ge&   i(\ga_k,2m_k+1) - i(\ga_k,1) - \frac{e(M_k)}{2} \nn\\
&=&  2N -  \frac{e(M_k)}{2} \nn\\
&\ge&  2N-n.    \nn\eaa
This proves  (\ref{eq4.55}).

  By the inequality (\ref{eq2.20}) in Theorem \ref{Thm2.2} as well as (\ref{eq4.54}),
\aa
i(\ga_k,2m_k) + \nu(\ga_k,2m_k)
&\le& i(\ga_k,2m_k+1) - i(\ga_k,1) + \frac{e(M_k)}{2}  \lb{eq4.67}\\
&=&  2N +  \frac{e(M_k)}{2}  \nn\\
&\le&  2N+n.  \nn\eaa
This proves  (\ref{eq4.56}) and the proof is complete. \hfill\qed\enddemo

  By Lemma \ref{Lem1.3}, the following consideration is useful for the study of
closed characteristics.

\specialnumber{4.4}\proclaim{Theorem}\lb{Thm4.4} Under the conditions of\ Theorem {\rm \ref{Thm4.3},}
further suppose
\ee  M_k=P_k^{-1}(N_1(1,1)\dm G_k)P_k   \lb{eq4.70}\eee
for some $P_k\in\Sp(2n)$ and $G_k\in\Sp(2n-2)$  holds for $1\le k\le q${\rm .}
Then {\rm (\ref{eq4.56})} can be improved to
\ee  i(\ga_k,2m_k) + \nu(\ga_k,2m_k) \le 2N+\frac{e(M_k)}{2}-1 \le 2N+n-1.\hskip.5in
    \lb{eq4.72}\eee
\et

\demo{Proof} In   Step 4 of the proof for  Theorem \ref{Thm4.3}, we use the inequality
(\ref{eq2.71}) in Theorem \ref{Thm2.3} as well as  (\ref{eq4.54}), and obtain
\aa
i(\ga_k,2m_k) + \nu(\ga_k,2m_k)
&\le& i(\ga_k,2m_k+1) - i(\ga_k,1) + \frac{e(M_k)}{2} - 1 \nn\\
&=&  2N +  \frac{e(M_k)}{2}-1  \nn\\
&\le&  2N+n -1.  \nn\eaa
This proves  (\ref{eq4.72}).   \enddemo
\pagebreak

\section{Proof of the main results}\lb{s5}

  Fix $\Sigma\in\H(2n)$ and $\alpha\in(1,2)$.

\specialnumber{5.1} \proclaim{Lemma}\lb{Lem5.1}  
 {\rm (i)} $\varrho_n(\Sigma)$ defined by {\rm (\ref{eq1.5})} does not depend on the choice
of $\alpha\in(1,2)${\rm .}
\bd\item{\rm (ii)} For any $s>0$ and $\Sigma\in\H(2n)${\rm ,} 
$\varrho_n(s\Sigma)=\varrho_n(\Sigma)${\rm .}\ebd
\el

\demo{Proof} (i) follows from Proposition I.7.5 of \cite{Ek3} and from a similar proof
for the index functions and splitting numbers.

\smallbreak
  (ii) Let $C$ be the convex compact set bounded by $\Sigma$. By definition,
\y
j_{sC}(x) &=& s^{-1}j_C(x),\\
H_{s\Sigma,\alpha}(x)&=&j_{sC}(x)^{\alpha}=s^{-\alpha}H_{\Sigma,\alpha}(x).
\ey
Let $(\tau,x)\in\J(\Sigma,\alpha)$. Then $\tau$ is the minimal period of $x$
and $H_{\Sigma,\alpha}(x(t))=1$. Set $y(t)=sx(s^{-2}t)$. Thus $s^2\tau$ is the minimal
period of $y$ and $H_{s\Sigma,\alpha}(y(t))=1$. Since $\dot x=JH^{\prime}_{\Sigma,\alpha}(x)$,
\y
\dot y(t)&=&s^{-1}JH^{\prime}_{\Sigma,\alpha}\lrp{x(s^{-2}t)}  \\
&=&s^{-\alpha}JH^{\prime}_{\Sigma,\alpha}\lrp{sx(s^{-2}t)}  \\
&=&JH^{\prime}_{s\Sigma,\alpha}(y(t)).\ey
So we have $(s^2\tau,y)\in\J(s\Sigma,\alpha)$. Therefore the map from $\J(\Sigma,\alpha)$
to $\J(s\Sigma,\alpha)$ defined by $(\tau,x)\mapsto (s^2\tau,sx(s^{-2}t))$ is a bijection.

  Let $\ga_x$ be the associated symplectic path of $(\tau,x)$. Then,
\y
\frac{d}{dt}\gamma_x(s^{-2}t)
&=& s^{-2}JH^{''}_{\Sigma,\alpha}\lrp{x(s^{-2}t))\gamma_x(s^{-2}t)}\\
&=& s^{-\alpha}JH^{''}_{\Sigma,\alpha}\lrp{sx(s^{-2}t))\gamma_x(s^{-2}t)}\\
&=& JH^{''}_{s\Sigma,\alpha}\lrp{y(t))\gamma_x(s^{-2}t)}.  \ey
So the associated symplectic path of $(s^2\tau,y)$ is $\ga_y(t)=\ga_x(s^{-2}t)$.
Hence the lemma follows. \enddemo

  Theorem \ref{Thm1.1} is contained in the following result.

\specialnumber{5.1}\proclaim{Theorem}\lb{Thm5.1} Let $\Sigma\in\H(2n)$ and $1<\alpha<2${\rm .}
Suppose $^{\#}\td{\J}(\Sigma)< +\infty${\rm .} Then 
\ee \min\lrc{^{\#}\V_{\infty}(\Sigma,\alpha),n}
    \ge \varrho_n(\Sigma)\ge \lrb{\frac{n}{2}}+1.         \lb{eq5.1}\eee
\et
\pagebreak

\demo{Proof}  The proof  is given by the following three claims.

\medbreak{\it Claim} 1. $^{\#}\V_{\infty}(\Sigma,\alpha)\ge\varrho_n(\Sigma)$.

 By Lemma \ref{Lem3.1}, there exist an integer $K\ge 0$ and
an injection map\break $p: \N+K \to \V_{\infty}(\Sigma,\alpha)\times\N$ such that
\ee
i(x,m) \le 2k-2+n \le i(x,m) + \nu(x,m) - 1    \lb{eq5.2}\eee
holds for any $k\in\N+K$, $(\tau,x)\in\J(\Sigma,\alpha)$
and $p(k)=([(\tau,x)],m)$. Denote the elements in
$\V_{\infty}(\Sigma,\alpha)$ by
$$  \V_{\infty}(\Sigma,\alpha)=\{[(\tau_j,x_j)]\mid j=1,\ldots,q\},  $$
where $(\tau_j,x_j)\in\J(\Sigma,\alpha)$ for $j=1,\ldots,q$. By Lemma
\ref{Lem1.3} and (\ref{eq3.4}) of Corollary \ref{Cor3.1}, we can apply
Theorems \ref{Thm4.3} and \ref{Thm4.4}, and obtain infinitely many
$(N,m_1,\ldots,m_q)\in \N^{q+1}$ such that
\aa \bigcap_{j=1}^q\G_{2m_j-1}(\tau_j,x_j)&\supset&[2N-\kappa_1,2N+\kappa_2],
   \lb{eq5.4}\\
\bigcup_{j=1}^q\I_{2m_j}(\tau_j,x_j)&\subset&[2N-n,2N-2+n],
   \lb{eq5.5}\eaa
where $\kappa_1=\kappa_1(\Sigma,\alpha)$ and $\kappa_2=\kappa_2(\Sigma,\alpha)$
are defined by (\ref{eq1.46}) and (\ref{eq1.47}) respectively. Recall that the $m^{\rm th}$
index interval $\I_m(\tau,x)$ and the $m^{\rm th}$ index jump $\G_m(\tau,x)$ of $(\tau,x)$
are defined by (\ref{eq1.31}) and (\ref{eq1.40}) respectively.  Since $\Sigma$ is strictly
convex, (\ref{eq1.24}) holds by Corollary  \ref{Cor1.2} for every $x_j$ with $1\le j\le q$.
Specifically, this yields
\ee  \kappa_2 \ge  n-1. \lb{eq5.6}\eee
Set $p(N-s+1)=([(\tau_{j(s)},x_{j(s)}], m(s))$ with $j(s)\in\{1,\ldots,q\}$  and
$m(s)\in\N$ for $s=1,\ldots,\varrho_n(\Sigma)$.
Then by definition of the map $p$ in Lemma \ref{Lem3.1}, 
\ee   i(x_{j(s)},m(s)) \le 2N - 2s + n  \le  i\lrp{x_{j(s)},m(s)}
 + \nu\lrp{x_{j(s)},m(s)} - 1. \quad \lb{eq5.6a}\eee
Because $\varrho_n(\Sigma)=[\frac{\kappa_1+n}{2}]$ and  (\ref{eq5.6}),
\aa
2N - \kappa_1
&\le&  2N - 2[\frac{\kappa_1+n}{2}] + n  \\
&\le&  2N - 2s + n  \\
&\le&  2N - 2 + n  \\
&<&  2N + \kappa_2,  \lb{eq5.6b}\eaa
for $s=1,\ldots,\varrho_n(\Sigma)$. From the definition (\ref{eq1.40}) \pagebreak  of the index jump,
(\ref{eq5.4}), (\ref{eq5.6a}), and (\ref{eq5.6b}), we obtain
\aa
&&\hskip-36pt i\lrp{x_{j(s)},2m_{j(s)}-1}+ \nu\lrp{x_{j(s)},2m_{j(s)}-1} - 1\lb{eq5.6c}\\
&&\hskip.25in < 2N - \kappa_1  \nn\\
&&\hskip.25in\le  2N - 2s + n \nn\\
&&\hskip.25in\le  i\lrp{x_{j(s)},m(s)} + \nu(x_{j(s)},m(s)) - 1,  \nn\eaa
and
\aa
i\lrp{x_{j(s)},m(s)}
&\le&  2N - 2s + n \lb{eq5.6d}\\
&<&  2N + \kappa_2   <  i\lrp{x_{j(s)},2m_{j(s)}+1}.  \nn\eaa
Comparing (\ref{eq5.6c}), (\ref{eq5.6d}), and (\ref{eq3.3}) of Corollary \ref{Cor3.1},
we obtain
\ee  2m_{j(s)} - 1 < m(s) < 2m_{j(s)} + 1.   \lb{eq5.6e}\eee
Hence $m(s)=2m_{j(s)}$;  i.e.,
\ee p(N-s+1)=\lrp{\lrb{(\tau_{j(s)},x_{j(s)}}},2m_{j(s)}),
         \qquad {\rm for}\; s=1,\ldots,\varrho_n(\Sigma). \hskip.25in \lb{eq5.7}\eee
Since the map $p$ is injective when $^{\#}\td{\J}(\Sigma)<+\infty$,
these $j(s)$'s are mutually different for $s=1,\ldots,\varrho_n(\Sigma)$. Therefore,
$$  q \ge \varrho_n(\Sigma).  $$
This proves Claim 1.
\vglue12pt

{\it Claim} 2. $\varrho_n(\Sigma)\le n$.

  In fact, if $\varrho_n(\Sigma)=[\frac{\kappa_1+n}{2}]>n$, we must have
$\kappa_1\ge n+2$. Denote by $([(\tau_j,x_j)],m)=p(N-n)$. Then 
$$  i(x_j,m) \le 2(N - n) - 2 + n \le i(x_j,m) + \nu(x_j,m) - 1. $$
By (\ref{eq5.4}),
\y
i(x_j,2m_j-1) + \nu(x_j,2m_j-1) - 1
&<& 2N - \kappa_1  \\
&\le& 2(N - n) - 2 + n  \\
&\le& i(x_j,m) + \nu(x_j,m) - 1,  \ey
and by (\ref{eq5.5}),
$$  i(x_j,m) \le 2(N - n) - 2 + n < 2N - n  \le i(x,2m_j).  $$
Hence by Corollary \ref{Cor3.1}, we have $2m_j-1<m<2m_j$.
This contradiction proves Claim 2.
\pagebreak

{\it Claim} 3.  $\varrho_n(\Sigma) \ge [\frac{n}{2}]+1$.

  Let $(\tau,x)\in\J(\Sigma,\alpha)$ and $\ga_x$ be its associated symplectic
path. In the following,  we estimate $i(x,1)+2S^+(x)-\nu(x,1)+n$.

  By Lemma \ref{Lem1.3}, (\ref{eq1.330a}) holds for some $P\in\Sp(2n)$ and
$M\in\Sp(2n-2)$. By $1^{\circ}$ and $2^{\circ}$ of Corollary 4.14 of
\cite{Lo7}, the splitting numbers are constant on each homotopy component
and symplectic additive. Thus by  (\ref{eq1.330a}),  
\ee  2S^+(x) - \nu(x,1) = 2S_{N_1(1,1)}^+(1) - \nu_1(N_1(1,1))
                    + 2S_{M}^+(1) - \nu_1(M).  \hskip.25in\lb{eq5.12}\eee
By  Theorem \ref{Thm6.6} below,  
$$  S_{N_1(1,a)}^+(1) = \left\{\matrix{1, &\quad {\rm if}\;\; a\ge 0, \cr
                                                              0, &\quad {\rm if}\;\; a< 0. \cr}\right. $$
Thus,
\ee  2S_{N_1(1,a)}^+(1) - \nu_1(N_1(1,a)) = a, \quad {\rm for}\;
       a=\pm 1, 0.   \lb{eq5.13}\eee
By Theorem \ref{Thm6.5},  
\ee  M \approx N_1(1,1)^{\dm p_-}\dm I_2^{\dm p_0}\dm N_1(1,-1)^{\dm p_+}\dm G,
\lb{eq5.14}\eee
for some nonnegative integers $p_-$, $p_0$, and $p_+$, and some symplectic
matrix $G$ satisfying $1\not\in \sigma(G)$. By (\ref{eq5.13}) and (\ref{eq5.14}),
we then obtain
\ee  2S_{M}^+(1) - \nu_1(M) \ge p_- -p_+ \ge - p_+\ge 1-n. \lb{eq5.15}\eee
From (\ref{eq1.24}), (\ref{eq5.12}), (\ref{eq5.13}) with $a=1$, and (\ref{eq5.15}), 
$$ i(x,1) + 2S^+(x) - \nu(x,1) + n  \ge  n + 1 + (1-n) + n = n + 2. $$
So Claim 3 holds and 
the proof  of Theorem \ref{Thm5.1} is complete. \enddemo

  Next we give the proof for the nondegenerate case.

\demo{Proof of Corollary {\rm \ref{Cor1.1}}}   By the definition of
$\varrho_n(\Sigma)$, (\ref{eq1.8}), and Theorem~\ref{Thm5.1}, we obtain
(\ref{eq1.9}).

  If $(\tau,x)\in\J(\Sigma,\alpha)$ is nondegenerate, we must have
$1\not\in\sg(M)$ in (\ref{eq5.12}). Thus by (\ref{eq5.12}) and (\ref{eq5.13}),
we obtain $2S^+(x) - \nu_{\tau}(x)=1$. This and   (\ref{eq1.24})
give  (\ref{eq1.8}). Now,  (\ref{eq1.9}) holds and by
(\ref{eq1.6}) this implies (\ref{eq1.10}). \enddemo

\specialnumber{5.1}\proclaim{{C}orollary}\lb{Cor5.1} Let $\Sigma\in\H(2n)$ and $\alpha\in (1,2)${\rm .}
Assume $^{\#}\td{\J}(\Sigma)< +\infty${\rm .} Then there exists an
element $[(\tau,x)]\in\V_{\infty}(\Sigma,\alpha)$ with $i_{\tau}(x)=n${\rm .}\ec

\demo{Proof} We use notation  introduced in the proof of Theorem \ref{Thm5.1}.
Let $([(\tau_j,x_j)],m)=p(N+1)$. By the definition of $p$, we obtain
$[(\tau_j,x_j)]\in\V_{\infty}(\Sigma,\alpha)$ and
\ee   i(x_j,m) \le 2N + n \le i(x_j,m) + \nu(x_j,m) - 1.  \lb{eq5.20}\eee
Thus by (\ref{eq5.4}) and (\ref{eq5.5}),
\aa
i(x_j,2m_j) + \nu(x_j,2m_j) - 1
&\le& 2N - 2 + n  \lb{eq5.21}\\
&<&  2N + n  \nn\\
&\le& 2N + \kappa_2 + 1  \nn\\
&\le& i(x_j,2m_j+1). \nn \eaa
Hence by Corollary \ref{Cor3.1}, comparing (\ref{eq5.20}) and
(\ref{eq5.21}) we obtain $2m_j<m\le 2m_j+1$.
Thus $m=2m_j+1$. By (\ref{eq5.20}) and (\ref{eq5.21}) again,
$$ i(x,2m_j+1) = 2N + n.  $$
Now by (\ref{eq4.54}) of Theorem \ref{Thm4.3},  
\medbreak \hfill $i(x,1) = i(x,2m_j+1) - 2N = n.  $ 
\enddemo

\demo{{R}emark {\rm5.1}} Although the global minimal point $(\tau,x)$ of
$f$ on $E$ satisfies $i_{\tau}(x)=n$, it is not clear whether
$[(\tau,x)]\in\V_{\infty}(\Sigma,\alpha)$.  \er

  Now we come to the proof of Theorem \ref{Thm1.2}. We will
prove the following stronger version. For every $(\tau,x)\in\J(\Sigma,\alpha)$,
the {\it elliptic height} $e(x)$ of $(\tau,x)$ is defined to be the elliptic height
$e(\ga_x(\tau))$, where $\ga_x$ is the associated symplectic path of $(\tau,x)$.

\specialnumber{5.2}\proclaim{Theorem}\lb{Thm5.2} Let $\Sigma\in\H(2n)$ and $\alpha\in(1,2)$. Suppose
$^{\#}\td{\J}(\Sigma)< +\infty$. By the notation introduced in
the proof of Theorem {\rm \ref{Thm5.1},} and specifically{\rm ,} the proof of Claim $1$ in the proof
of  Theorem {\rm \ref{Thm5.1},} for each $s=1,\ldots,\varrho_n(\Sigma)${\rm ,}
there exists a unique $j(s)\in\{1,\ldots,q\}$ such that
$p(N-s+1)=([(\tau_{j(s)},x_{j(s)})],2m_{j(s)})$.   Then
\ee  e(x_{j(s)}) \ge 2|n-2s+1|+2, \qquad
      \hbox{for all }s=1,\ldots,\varrho_n(\Sigma).\hskip.5in \lb{eq5.22}\eee
Specially{\rm ,} $[(\tau_{j(1)},x_{j(1)})]$ is an elliptic element in
$\V_{\infty}(\Sigma,\alpha)${\rm .}  \et

\demo{Proof} By the definition of the injection map $p$ in Lemma \ref{Lem3.1},
we have
\ee
i(x_{j(s)},2m_{j(s)}) \le 2N - 2s + n
  \le  i(x_{j(s)},2m_{j(s)}) + \nu(x_{j(s)},2m_{j(s)}) - 1.  \lb{eq5.23}\eee
By (\ref{eq2.71}) of  Theorem \ref{Thm2.3} and (\ref{eq4.54}) in
Theorem \ref{Thm4.3},  
\aa
\frac{e(\tau_{j(s)},x_{j(s)})}{2}
&\ge& i\lrp{x_{j(s)},2m_{j(s)}} + \nu\lrp{x_{j(s)},2m_{j(s)}}\lb{eq5.24}\\
&&\ - i\lrp{x_{j(s)},2m_{j(s)}+1}
         + i\lrp{x_{j(s)},1} + 1  \nonumber\\
&=& i\lrp{x_{j(s)},2m_{j(s)}} + \nu\lrp{x_{j(s)},2m_{j(s)}} - 2N + 1  \nn\\
&\ge& 2N - 2s + n +1 - 2N + 1  \nn\\
&=& n - 2s + 2,  \lb{eq5.23a}\eaa
where (\ref{eq5.23a}) follows from the right-hand   inequality in (\ref{eq5.23}).

  On the other hand, by (\ref{eq2.72}) in Theorem \ref{Thm2.3},  
\aa
\frac{e(\tau_{j(s)},x_{j(s)})}{2}
&\ge& i\lrp{x_{j(s)},2m_{j(s)}+1} + \nu\lrp{x_{j(s)},2m_{j(s)}+1}   \nn\\
 & &    -\ i\lrp{x_{j(s)},2m_{j(s)}} - i\lrp{x_{j(s)},1} - \nu\lrp{x_{j(s)},1}. \nn\eaa
By (\ref{eq4.52}) and (\ref{eq4.54}) of Theorem \ref{Thm4.3} and the left-hand  
inequality in (\ref{eq5.23}), 
\aa
\frac{e(\tau_{j(s)},x_{j(s)})}{2}
&\ge& 2N - (2N - 2s + n)    \\
&=& 2s - n. \lb{eq5.25}\eaa
Now (5.24) and (\ref{eq5.25}) yield (\ref{eq5.22}).  \enddemo

 The following theorem implies Theorem \ref{Thm1.3}. Note that by Theorem
\ref{Thm5.1}  $\varrho_n(\Sg)-1\ge [\frac{n}{2}]$ always.

\specialnumber{5.3}\proclaim{Theorem}\lb{Thm5.3} Let $\Sigma\in\H(2n)$ with $n\ge 2$ and
$\alpha\in(1,2)${\rm .}
Suppose $^{\#}\td{\J}(\Sigma)< +\infty${\rm .} Then there exist at least
$\varrho_n(\Sg)-1$ elements in $\V_{\infty}(\Sigma,\alpha)$ such that
each such an element $[(\tau,x)]$ satisfies
\ee  \hat{i}(x,1) \in \R\bs\Q.      \lb{eq5.26}\eee
Specifically{\rm ,} such a closed characteristic $[(\tau,x)]$ must possess at least a
rotator on $\U$ in the sense of Theorem {\rm \ref{Thm2.1}} with an irrational multiple
of $2\pi$ rotation angle{\rm .}
\et

\demo{Proof} Since $n\ge 2$, we have $\varrho_n(\Sg)-1\ge 1$. Let
$$  \V_{\infty}(\Sigma,\alpha)=\lrc{\lrb{(\tau_j,x_j)}\,|\,1\le j\le q }.  $$
By (\ref{eq4.58}) and Corollary \ref{Cor2.1}, the $D_j$ defined by (\ref{eq4.1})
according to Step 1 of the proof of Theorem  \ref{Thm4.3} satisfies
\ee  D_j=\hat{i}(x_j,1), \qquad \hbox{for all }j=1,\ldots,q.  \lb{eq5.28}\eee
By Step 1 of the proof of Theorem \ref{Thm5.1}, we obtain
$(N, m_1, \ldots, m_q)\in\N^{q+1}$ such that the injection map $p$ of
Lemma \ref{Lem3.1} satisfies  (\ref{eq5.7}). Reordering elements in
$\V_{\infty}(\Sigma,\alpha)$ to simplify notation, we can assume
\ee  p(N-j+1)=([\tau_j, x_j],2m_j), \qquad \hbox{for all }j=1,\ldots,\varrho_n(\Sg).\hskip.25in
      \lb{eq5.30}\eee
Thus by Lemma \ref{Lem3.1} and (\ref{eq1.19a}),  
\aa
2m_kD_k - 2m_jD_j
&=& 2m_k\hat{i}(x_k,1) - 2m_j\hat{i}(x_j,1)  \lb{eq5.32}\\
&=& \hat{i}(x_k,2m_k) - \hat{i}(x_j,2m_j)  \nn\\
&<& 0, \qquad  {\rm if}\;\; 1\le j<k\le \varrho_n(\Sg).\nn\eaa

  Now it suffices to prove that among the first $\varrho_n(\Sg)$ of
$D_j$'s, at most one of them is rational. We prove this claim indirectly
by assuming
\ee   D_j \;\;\; {\rm and}\;\;\; D_k \in \Q, \lb{eq5.33}\eee
for some $j$ and $k$ satisfying $1\le j<k\le \varrho_n(\Sg)$. Then in our
choice of  $(N,m_1,\ldots,m_q)$ in the proof of Theorem \ref{Thm4.1}, we
require that $N\in \N$ further satisfies
\ee   \frac{N}{MD_i}\in\N, \quad {\rm for}\;\; i=j,\;{\rm and}\; k.
    \lb{eq5.34}\eee
Then the closure of the set
$\{\{Nv\}\,|\,N\in\N\;{\rm satisfies}\;(\ref{eq5.34})\}$
is still a closed additive subgroup of  $\T^h$ for some $h\in\N$. Thus
(\ref{eq4.38}) holds in the $j^{\rm th}$ and $k^{\rm th}$ coordinates for infinitely
many $N$, if we choose $\chi$ to be a vertex of $[0,1]^h$ with
$\chi_j=\chi_k=0$. This implies
$$ \left|\lrc{\frac{N}{MD_i}}\right| = 0, \qquad {\rm for}\;\; i=j, k.  $$
Therefore by the definition (\ref{eq4.14}) of $m_i$'s,  
$$  m_i = \lrb{\frac{N}{MD_i}}M=\frac{N}{D_i}, \qquad {\rm for}\;\; i=j, k.  $$
This yields
\ee  m_jD_j=N =m_kD_k, \lb{eq5.35}\eee
  contradicts   (\ref{eq5.32}), and completes the proof of the theorem. \enddemo

  The following theorem shows the existence of  multiple elliptic orbits
and implies Theorem \ref{Thm1.4}.

\specialnumber{5.4}\proclaim{Theorem}\lb{Thm5.4} Let $\Sigma\in\H(2n)$ with $n\ge 2$ and
$\alpha\in(1,2)${\rm .}
Suppose $^{\#}\tilde\J(\Sigma)< +\infty$ and
\ee   ^{\#}\V_{\infty}(\Sigma,\alpha)\le 2\varrho_n(\Sigma)-2.   \lb{eq5.36}\eee
Then there exist at least two elliptic elements in $\V_{\infty}(\Sigma,\alpha)${\rm .}

  In particular{\rm ,}  by {\rm (\ref{eq1.7}),} there are at least two elliptic elements in
$\V_{\infty}(\Sigma,\alpha)${\rm ,} provided
\ee   ^{\#}\V(\Sigma,\alpha)\le 2\lrb{\frac{n}{2}}.  \lb{eq5.38}\eee
\et

\demo{Proof}  Let $\V_{\infty}(\Sigma,\alpha)=\{[(\tau_j,x_j)]\,|\,1\le j\le q\}$. Let
$v\in\R^k$ be the vector given by (\ref{eq4.36}) according to quantities in
(\ref{eq4.58a})--(\ref{eq4.58b}) in the proof of Theorem~\ref{Thm4.3}. Let $A(v)$
be the set given by Theorem \ref{Thm4.2} according to this $v$.

  We use notation  introduced in the proof of Theorem \ref{Thm5.3}. Note that by
Theorem \ref{Thm5.1} and the fact that $n\ge 2$,  $q\ge \varrho_n(\Sg) \ge 2$. By
Theorem~\ref{Thm5.3}, $D_j\in \R\bs\Q$ holds for at least $[n/2]$ integers in
$\{1,\ldots,q\}$. Therefore by Theorem~\ref{Thm4.2},  
\ee  \dim A(v) \ge 1, \qquad 0\not\in A(v).  \lb{eq5.40}\eee

  We continue our proof in two steps.

\medbreak {\it  Step} 1.  For a given $a\in A(v)$, we define
$\chi\equiv\chi(a)=(\psi(a_1),\ldots,\psi(a_k))$ by (\ref{eq4.44}). Let
$(N, m_1, \ldots, m_q)\in\N^{q+1}$ be given in Step 1 of  the proof  for Theorem
\ref{Thm5.1} via an application of Theorems \ref{Thm4.1} to  \ref{Thm4.3} with
this $\chi=\chi(a)$ used in their proofs. By the proof of Theorem \ref{Thm5.1},
there exists a subset  $L(a,N)$ of  $\{1,\ldots,q\}$ such that
\bd
\item{(i)} $1\in L(a,N)$ and $^{\#}L(a,N)= \varrho_n(\Sg)$,
\item{(ii)} $p(N)=([(\tau_1,x_1),m_1)]$ and $x_1$ is elliptic;
\item{(iii)} For each $s=2,\ldots, \varrho_n(\Sg)$, there is a
$j(s)\in L(a,N)\bs\{1\}$ such that $$p(N-s+1)=([\tau_{j(s)},x_{j(s)}],m_{j(s)}).$$
\ebd

  By Lemma \ref{Lem3.1}, we have
\ee  m_1D_1-m_jD_j >0, \quad \hbox{for all }j\in L(a,N)\bs\{1\}.  \lb{eq5.42}\eee

Assuming that there exists only one elliptic element in
$\V_{\infty}(\Sigma,\alpha)$, we always get the same elliptic element
$[(\tau_1,x_1)]\in\V_{\infty}(\Sigma,\alpha)$ for all $a\in A(v)$.

  By (\ref{eq5.40}), we have $-a\in A(v)$. By   definition (\ref{eq4.44}),
  we have   $\chi(-a)\not=\chi(a)$. Thus by the proof of Theorem \ref{Thm5.1},
for this vertex $\chi(-a)$ we get $(\td{N},\td{m_1},\ldots,\td{m_q})\break\in\N^{q+1}$ and another subset
$L(-a,\td{N})$ of $\{1,\ldots,q\}$ such that Theorems \ref{Thm5.1}--\ref{Thm5.3} hold, specially the above
(i)--(iii) still hold
correspondingly.

  By the assumption $q\le 2\varrho_n(\Sigma)-2$ of (\ref{eq5.36}), 
\ee  1 \le \; ^{\#}((L(a,N)\cap L(-a,\td{N}))\bs\{1\})\le \varrho_n(\Sg).
    \lb{eq5.44}\eee

  Let
\ee  t_0 = \frac{\delta_1}{6(|a|+1)(M\Lambda+1)}, \qquad
                \Lambda = \max_{1\le j\le q}D_j.  \lb{eq5.50}\eee
Note that in the proof of  Theorem \ref{Thm4.1} we can further require $N\in \N$ such
that the vector $\{Nv\}-\chi(a)$ defined by (\ref{eq4.38}) and  (\ref{eq4.36}) are located in a sufficiently
small neighborhood inside the open ball in $V$ with radius
$\delta_1/(6M\Lambda+1)$ and centered at $at_0$; i.e.,
\aa
  && \{Nv\} - \chi(a) \in V,  \lb{eq5.51}\\
   && |\{Nv\} - \chi(a) - at_0| < \frac{\delta_1}{6M\Lambda+1},
  \lb{eq5.52}\eaa
where $V=T_0\pi^{-1}(\ol{\{\{Nv\}\,|\,N\in\N\}})$ is as defined in Theorem
\ref{Thm4.2}. Under this requirement, we still have (\ref{eq4.38}):
\aa
  \left|\lrc{\frac{N}{MD_i} }-\chi_i\right| &<& |a_it_0| + \frac{\delta_1}{6M\Lambda+1}
             \le \frac{\delta_1}{3M\Lambda},  \lb{eq5.54}\\
  \left|\lrc{\frac{N\alpha_{i,j}}{D_i} }-\chi_{i,j}\right| &<&  |a_{i,j}t_0|
        + \frac{\delta_1}{6M\Lambda+1} \le \frac{\delta_1}{3},  \lb{eq5.56}\eaa
for infinitely many $N\in \N$. Here $\chi_i=\psi(a_i)$ and
$\chi_{i,j}=\psi(a_{i,j})$; the function $\psi$ is defined in Theorem \ref{Thm4.2}.

\vglue12pt {\it Claim}.  $a_1D_1-a_jD_j=0$ for all
$j\in (L(a,N)\cap L(-a,\td{N}))\bs\{1\}$.  
  
Assuming the claim does not hold, we prove it by contradiction. In fact, we
can further require $N\in \N$ so that the following also holds:
$$  \left|\lrc{\frac{N}{MD_k} } - \chi_k - a_kt_0\right|
  < \frac{t_0}{3\Lambda}\min_{a_iD_i-a_jD_j\ne 0}|a_iD_i-a_jD_j|,\qquad
                \hbox{for all }1\le k\le q. $$
By our choice of $(N,m_1, \ldots,m_q)$, we have
\aa
&&\lb{eq5.58}\\
m_1D_1-m_jD_j
&=& M\lrp{\lrb{\frac{N}{MD_1}}+\chi_1}D_1-M\lrp{\lrb{\frac{N}{MD_j}}+\chi_j}D_j \nn\\
&=& M\left(\lrp{\chi_1-\lrc{\frac{N}{MD_1} }}D_1
        - \lrp{\chi_j-\lrc{\frac{N}{MD_j} }}D_j\right)   \nn\\
&=& - Mt_0(a_1D_1-a_jD_j)
          + M\left(\lrp{\chi_1-\lrc{\frac{N}{MD_1} }+a_1t_0}D_1\right.      \lb{eq5.60}\\
& & \qquad  \left.-\lrp{\chi_j-\lrc{\frac{N}{MD_j} }+a_jt_0}D_j\right).
  \nn\eaa

  For any $j\in (L(a,N)\cap L(-a,\td{N}))\bs\{1\}$, by (\ref{eq5.42}) and
(\ref{eq5.60}),
\ee a_1D_1\le a_jD_j.   \lb{eq5.62}\eee

  For any $j>1$ satisfying $j\in L(a,N)\cap L(-a,\td{N})$, we apply the above
argument to $a$ and $-a$ respectively and obtain $a_1D_1\le a_jD_j$ and
$-a_1D_1\le -a_jD_j$. Therefore $a_1D_1=a_jD_j$, and the claim is proved.
\vglue12pt

{\it Step} 2.  Set
\aa
 V_j &=& \{a\in V\mid a_1D_1=a_jD_j\},\quad j=2,\ldots,q,  \lb{eq5.63}\\
 B(v) &=& A(v)\bs\bigcup_{V_j\ne V,\;j=2,\ldots,q}V_j.  \lb{eq5.64}\eaa

  Since $\dim A(v)\ge 1$, by Theorem \ref{Thm4.2}, $A(v)$ is obtained from $V$ by
deleting finitely many proper linear subspaces of $V$, and so is $B(v)$. Hence $B(v)$ is nonempty.

  Now we choose an $a\in B(v)$ and $j\in L(a,N)\cap L(-a,\tilde N)\bs\{1\}$.
By the Claim in   Step 1, we have $a_1D_1=a_jD_j$. By   definitions of
$a\in B(v)$ and $j$ we have $V_j=V$.

  By  (\ref{eq5.51}), the vector $\{Nv\}-\chi$ defined by (\ref{eq4.38}) belongs to
$V$, and thus belongs to $V_j$. Then by   definition (\ref{eq5.63}) of $V_j$, this
implies
\ee   0 = (\{Nv_1\}-\chi_1)D_1-(\{Nv_j\}-\chi_j)D_j.  \lb{eq5.65}\eee
By   definition (\ref{eq4.36}) of $v$,   
$$ 0 = \lrp{\lrc{\frac{N}{MD_1} }-\chi_1}D_1-\lrp{\lrc{\frac{N}{MD_j} }-\chi_j}D_j.  $$
By (\ref{eq5.58}), this implies $m_1D_1=m_jD_j$. Then it contradicts  
(\ref{eq5.42}) and completes the proof.  \enddemo

\vglue-12pt
\section{Appendix 1. The Maslov-type index and its iteration theory}\lb{s6}
\vglue-8pt

In this section we give a brief review of the Maslov-type index and its iteration
theory for symplectic matrix paths. We use notation  already introduced in  
Section \ref{s1}.

\demo{Definition {\rm 6.1 (cf.\ \cite{Lo1} and \cite{Lo7})}} For $\tau>0$ and $\om\in\U$, given
two paths $\ga_0$ and $\ga_1\in\P_{\tau}(2n)$, if there exists a map
$\dl\in C([0,1]\times [0,\tau],\Sp(2n))$ such that
$\dl(0,\cdot)=\ga_0(\cdot)$, $\dl(1,\cdot)=\ga_1(\cdot)$, $\dl(s,0)=I$,
and $\nu_{\om}(\dl(s,.))$ is constant for $0\leq s\leq 1$, then $\ga_0$
and $\ga_1$ are {\it $\om$-homotopic on $[0,\tau]$ along $\dl(\cdot,\tau)$ }
and  $\ga_0\sim_{\om}\ga_1$. If $\ga_0\sim_{\om}\ga_1$ on $[0,\tau]$
along $\dl(\cdot,\tau)$ for all $\om\in\U$, then $\ga_0$ and $\ga_1$ are
{\it homotopic on $[0,\tau]$ along $\dl(\cdot,\tau)$ } and   we write
$\ga_0\sim\ga_1$. \enddemo

  We define $2\times 2$ and $4\times 4$ matrices
\aa
 D(a) &=& \left(\matrix{a & 0 \cr
                                     0 & \frac{1}{a} \cr}\right), \qquad
   R(\th)=\left(\matrix{\cos\th & -\sin\th \cr
                                   \sin\th &  \cos\th \cr}\right),  \nn\\
  N_1(\lm,c)&=&\left(\matrix{\lm & c\cr
                              0 & \lm\cr}\right), \qquad
   N_2(\om,b)=\left(\matrix{R(\th) & b\cr
                         0        & R(\th)\cr}\right),  \nn\eaa
for $a, \lm, c, \th\in\R$, $\om\in\U$, and $b\in\Li(\R^2)$. They are called
{\it basic normal forms} of symplectic matrices in \cite{Lo7} when
$a=\pm 2$, $\th\in\R$, $\lm=\pm 1$, $c\in\R$,
$\om=e^{\th\sqrt{-1}}\in \U\bs\R$ with $\th\in\R$;
$b=\left(\matrix{b_1 & b_2\cr
                           b_3 &  b_4\cr}\right)$ is a $2\times 2$ real
matrix with $b_2-b_3\not=0$. In \cite{Lo7}, a basic normal form matrix
$M\in\Sp(2n)$ is called {\it  trivial} if
$\sg(MR((t-1)\ep)^{\dm n})\cap\U=\emptyset$ for $\ep>0$ small
enough and $t\in [0,1)$. Note that the normal forms $N_1(1,-1)$,
$N_1(-1,1)$, $D(\pm 2)$, $N_2(\om,b)$ and $N_2(\ol{\om},b)\in \M_{\om}^1(4)$
with $\om=\exp(\th\sqrt{-1})\in\U\bs\R$ and $(b_2-b_3)\sin\th>0$ are trivial,
and all other basic normal form matrices are {\it nontrivial}. Here using
Definition \ref{Def1.2} we define
$$ \M_{\om}^1(2n) = \{M\in\Sp(2n)\,|\, \nu_{\om}(M)=1\}. $$

  For any $\tau>0$, $n\in\N$, and $k\in\Z$, following \cite{Lo1}, we define
\aa
  \phi_{\th,\tau}(t) &=& R\lrp{\frac{\th t}{\tau}}, \quad \hbox{for all }\;t\in [0,\tau], \nn\\
  \chi_{\tau}(t) &=& D\lrp{1+\frac{t}{\tau}}, \quad \hbox{for all }\;t\in [0,\tau]; \nn\eaa
We also denote by $N_1(e^{\theta\sqrt{-1}},b)=R(\theta)$ for any $b\in \R$ and $\theta\in \R\setminus
(\pi \Z)$. 

The following theorem gives a characterization of our index \hbox{function theory.}

 \proclaimtitle{Theorem 2.11 of  \cite{Lo7}}
\specialnumber{6.1}\proclaim{Theorem}\lb{Thm6.1}  For any
$\tau>0$ and
$\om=\exp(\th\sqrt{-1})\break\in\U${\rm ,} there exists a unique function
$i_{\om}:\cup_{n\in\N}\P_{\tau}(2n)\to\Z$ satisfying the following
five axioms\/{\rm :}\/
\medbreak
  $1^{\circ}\,$ {\rm (}\/Homotopy invariant\/{\rm ).}\/  For $\ga_0$ and
$\ga_1\in\P_{\tau}(2n)$, if $\ga_0\sim_{\om}\ga_1$ on $[0,\tau]${\rm ,} then
$$ i_{\om}(\ga_0) = i_{\om}(\ga_1). $$

\medbreak

  $2^{\circ}\,$ {\rm (}Symplectic additivity\/{\rm ).}   For any
$\ga_i\in\P_{\tau}(2n_i)$ with $i=0$ and $1${\rm , }
$$ i_{\om}(\ga_0\dm\ga_1) = i_{\om}(\ga_0) + i_{\om}(\ga_1). $$

\medbreak
  $3^{\circ}\,$ {\rm (}Clockwise continuity\/{\rm ).}  For any $\ga\in\P_{\tau}(2)$ and $\om\in\U$
satisfying $\ga(\tau)=N_1(\om,b)${\rm ,} there exists a $\th_0>0$ such that
$$ i_{\om}((\ga(\tau)\phi_{-\th,\tau})\ast\ga)
     = i_{\om}(\ga), \qquad \hbox{for all }\;0<\th\leq \th_0.   $$

\medbreak
  $4^{\circ}\,$ {\rm (}Counterclockwise jumping\/{\rm ).}  For any $\ga\in\P_{\tau}(2)$ and
$\om\in\U$ satisfying $\ga(\tau)=N_1(\om,b)$,  there exists a $\th_0>0$ such that
$$ i_{\om}((\ga(\tau)\phi_{\th,\tau})\ast\ga)
  = i_{\om}(\ga)+1, \qquad \hbox{for all }\;0<\th\leq \th_0.   $$

\medbreak
  $5^{\circ}\,$ {\rm (}Normality\/{\rm ).} 
$$ i_{\om}(\chi_{\tau}) = 0.   $$
\et

\proclaimtitle{Theorem 2.13 of \cite{Lo7}}
\specialnumber{6.2 }\proclaim{Theorem}\lb{Thm6.2}   For any two
paths
$\ga_0$ and\break $\ga_1\in \P_{\tau}(2n)$ with
$i_{\om}(\ga_0)=i_{\om}(\ga_1)${\rm ,} suppose that there exists a
continuous\break  path $h:[0,1]\to\Sp(2n)$ such that $h(0)=\ga_0(\tau)${\rm ,}
$h(1)=\ga_1(\tau)${\rm ,} and\break
$\dim_{\C}\ker_{\C}(h(s)-\om I)=\nu_{\tau,\om}(\ga_0)$
for all $s\in [0,1]${\rm .} Then $\ga_0\sim_{\om}\ga_1$ on $[0,\tau]$
along $h${\rm .} \et

\proclaimtitle{Theorem 1.4 of \cite{Lo7}}
\specialnumber{6.3}\proclaim{Theorem}\lb{Thm6.3}  For any
$\tau>0${\rm ,}
$\ga\in\P_{\tau}(2n)$, and $k\in\N${\rm ,}  
\ee
  i(\ga,k) = \sum_{\om^k=1}i_{\om}(\ga),  \quad
  \nu(\ga,k) = \sum_{\om^k=1}\nu_{\om}(\ga).  \nn\eee
\et

\proclaimtitle{Theorem 1.5 of  \cite{Lo7}}
\specialnumber{6.4}\proclaim{Theorem}\lb{Thm6.4}   For any
$\tau>0$ and
$\ga\in\P_{\tau}(2n)${\rm ,}
\ee
  \hat{i}(\ga,1) \equiv \lim_{k\to +\infty}\frac{i(\ga,k)}{k}
         = {1\over{2\pi}}\int_{\U}i_{\om}(\ga)d\om.  \nn\eee
Specially{\rm ,} $\hat{i}(\ga,1)$ is always a finite real number{\rm ,} which is called the
{\rm mean index} per $\tau$ of $\ga$. \et

\demo{Definition {\rm 6.2 (Definition 1.1 of \cite{Lo7})}} For any
$M\in\Sp(2n)$, define the {\it homotopy set} of $M$ in $\Sp(2n)$ by
\aa
  \Om(M)&\hskip-7pt=\hskip-7pt&\{N\in\Sp(2n) \,|\, \sg(N)\cap\U=\sg(M)\cap\U,\,\mbox{and}\; \nn\\
   &\hskip-7pt\hskip-7pt&  \dim_{\C}\ker_{\C}(N-\lm I) = \dim_{\C}\ker_{\C}(M-\lm I),\;
      \hbox{for all }\,\lm\in\sg(M)\cap\U\}.  \nn\eaa
We denote by $\Om^0(M)$ the path connected component of $\Om(M)$ which
contains $M$, and call it the {\it homotopy component} of $M$ in $\Sp(2n)$.
Also, $N\approx M$ if $N\in \Om^0(M)$. This is an equivalent relation
in $\Sp(2n)$. The equivalent class of $M$ in $\Sp(2n)$ is called the
{\it homotopy type} of $M$. \enddemo

\proclaimtitle{Theorem 7.8 of \cite{Lo7}}
\specialnumber{6.5}\proclaim{Theorem}\lb{Thm6.5}  For any
$M\in\Sp(2n)${\rm ,}   we have
\aa  
\noalign{\vskip-6pt}
&&\lb{eq6.10}\\M
&\approx&N_1(1,1)^{\dm p_-}\dm I_{2p_0}\dm N_1(1,-1)^{\dm p_+}\dm N_1(-1,1)^{\dm q_-}
              \dm (-I_{2q_0})\nn\\&&\; \dm N_1(-1,-1)^{\dm q_+}   
\dm R(\th_1)\dm\cdots\dm R(\th_r)\dm N_2(\om_1,u_1)\dm\cdots
        \dm N_2(\om_{r_{\ast}},u_{r_{\ast}}) \nn\\
&&\;\dm N_2(\lm_1,v_1)\dm\cdots\dm N_2(\lm_{r_0},v_{r_0})\dm M_k,
       \nn\eaa
where $p_-, p_0, p_+, q_-, q_0, q_+, r, r_{\ast}$, $r_0${\rm ,} and $k$ are nonnegative
integers\/{\rm ;} $\om_j=e^{\sqrt{-1}\alpha_j}${\rm ,} $\lm_j=e^{\sqrt{-1}\beta_j}$\/{\rm ;} the
real numbers $\th_j$ with $1\leq j\leq r${\rm ,} $\alpha_h$ with
$1\leq h\leq r_{\ast}${\rm ,} $\beta_l$ with $1\leq l\leq r_0${\rm ,} are in
$(0,\pi)\cup (\pi,2\pi)$ provided the corresponding $r${\rm ,} $r_{\ast}${\rm ,} or
$r_0>0$ respectively\/{\rm ;} $N_2(\om_j,u_j)$\/{\rm '}\/s are nontrivial and
$N_2(\lm_j,v_j)$\/{\rm '}\/s are trivial basic normal forms\/{\rm ;} $M_k=D(2)^{\dm k}$ or
$D(-2)\dm D(2)^{\dm (k-1)}${\rm .} All these integers{\rm ,} real numbers{\rm ,} and basic normal
form matrices are uniquely determined by $M${\rm .} It holds that 
\ee  p_-+p_0+p_++q_-+q_0+q_-+r+2r_{\ast}+2r_0+k=n. \lb{eq6.11}\eee
\et

{\it Definition {\rm 6.3  (Definition 4.8 of \cite{Lo7})}}.
For any basic normal form $M\in\Sp(2n)$ and
$\om\in\U\,\cap\,\sg(M)$, we define the
{\it ultimate type} $(p,q)$ of $\om$ for $M$ to be its usual Krein type
if $M$ is nontrivial, and to be $(0,0)$ if $M$ is trivial. When
$\om\in\U\bs\sg(M)$ with $M\in\Sp(2n)$, we define the ultimate type
of $\om$ for $M$ to be $(0,0)$. For any $M\in\Sp(2n)$, by Theorem
 6.5 there exists a $\dm$-product expansion (6.3) in the homotopy component $\Om^0(M)$ of $M$ as:
$$
M\approx M_1\dm M_2\dm\cdots \dm M_k\dm M_0,
$$
 where
each
$M_i$ is a basic normal form for
$1\leq i\leq k$ and
$\sg(M_0)\cap\U=\emptyset$. Denote the ultimate type of $\om$ for $M_i$
by $(p_i,q_i)$ for $0\leq i\leq k$. Let $p=\sum_{i=0}^kp_i$ and
$q=\sum_{i=0}^kq_i$. We define the {\it ultimate type} of $\om$ for $M$
by $(p,q)$. 
\vglue3pt

  The following theorem characterizes the splitting numbers algebraically.

\proclaimtitle{Lemma 4.5 and Theorem 4.11 of
\cite{Lo7}}
\specialnumber{6.6}\proclaim{Theorem}\lb{Thm6.6}  \hglue-7pt For any
$\om\in\U$~and $M\in\Sp(2n)${\rm ,} both $S_N^+(\om)$ and $S_N^+(\om)$ are
constants for any $N\in\Om^0(M)${\rm ,} and  
\ee  S_M^+(\om) = p \quad {\it and}\quad S_M^-(\om) = q, \nn\eee
where $(p,q)$ is the ultimate type of $\om$ for $M${\rm .} \et
\pagebreak

\proclaimtitle{Theorem 1.3 of \cite{Lo9}}
\specialnumber{6.7}\proclaim{Theorem}\lb{Thm6.7}   For
$\tau>0${\rm ,} let
$\ga\in
\P_{\tau}(2n)${\rm .} In
Theorem {\rm \ref{Thm6.5}} we let $M=\ga(\tau)$ and use the  notation   there{\rm .} Then for any
$m\in\N$,
\aa  \qquad i(\ga,m)
  &\hskip-7pt=\hskip-7pt& m(i(\ga,1)+p_-+p_0-r ) + 2\sum_{j=1}^rE\lrp{{{m\th_j}\over{2\pi}}} - p_-
\lb{eq6.13}\\
      &\hskip-7pt\hskip-7pt &    - p_0 - {{1+(-1)^m}\over 2}(q_0+q_+)   
  - r + 2\lrp{\sum_{j=1}^{r_{\ast}}\phi\lrp{\frac{m\alpha_j}{2\pi}} - r_{\ast}},  \nn \\
  \nu(\ga,m)
  &\hskip-7pt=\hskip-7pt& \nu(\ga,1) + {{1+(-1)^m}\over 2}(q_-+2q_0+q_+)
                 + 2\vf(m,\ga(\tau)),  \lb{eq6.14}\eaa
where 
\begin{eqnarray}   \vf(m,\ga(\tau))& = &\lrp{r-\sum_{j=1}^r\phi\lrp{\frac{m\th_j}{2\pi}}}
         + \lrp{r_{\ast}-\sum_{j=1}^{r_{\ast}}\phi\lrp{\frac{m\alpha_j}{2\pi}}}
 \lb{eq6.15} \\
&&            + \lrp{r_0-\sum_{j=1}^{r_0}\phi\lrp{\frac{m\beta_j}{2\pi}}}.   \nn\end{eqnarray}
Here the functions $E(\cdot)$ and $\phi(\cdot)$ are as defined at the beginning of
Section {\rm \ref{s2}.}
\et

\setcounter{equation}{0}

\section{Appendix 2. The Fadell-Rabinowitz $S^1$-cohomological index.\\
(By John Mather)}\lb{s7}

  For a principal $U(1)$-bundle $E\to B$, the Fadell-Rabinowitz index of $E$
is defined to be
$$  \sup\{k\;|\, c_1(E)^{k-1}\not= 0\}, $$
where $c_1(E)\in H^2(B,\Q)$ is the first rational Chern class. For a 
$U(1)$-space, i.e., a topological space $X$ with a $U(1)$-action, the
Fadell-Rabinowitz index is defined to be the index of the bundle
$$  X\times S^{\infty} \to  X\times_{U(1)}S^{\infty},  $$
where $S^{\infty}\to CP^{\infty}$ is the universal $U(1)$-bundle.

\vglue4pt {\it Acknowledgements}. The authors sincerely thank Professors Helmut
Hofer, John Mather, and the referee  for their careful reading and many
valuable comments on this paper, which make it more readable, with special  thanks to
  Professor John Mather for  Appendix 2,
and to the referee  for pointing out the equivalent definition (\ref{eq1.17})
of the index function in the degenerate case.
\pagebreak

\bye